{
\documentclass[doublespacing]{elsart}
\usepackage{amsmath,amssymb,graphicx,mathdots,natbib}
\bibpunct{[}{]}{,}{n}{,}{,}				

\newcommand{\eq}[1]{\eqref{#1}}
\newcommand{\fig}[1]{Fig.\ \ref{#1}}
\newcommand{\gaspar}{GASpAR}
\newcommand{\Rn}{{\sf R\mspace{-2mu}e}}
\newcommand{\sr}[1]{Section \ref{#1}}
\def\dof{d.o.f.} 
\def\dofs{d.o.f.} 
\def\cf{cf.} 
\def\eg{e.g.}  
\def\nse{Navier-Stokes equations\ }

\newcommand{\be}{\begin{equation}}
\newcommand{\ee}{\end{equation}}
\def\u{{\vec u}} \def\v{{\vec v}} \def\x{{\vec x}}

\def \pmbtext#1{\leavevmode \setbox0\hbox{#1}
     \kern-0,2pt \copy0 \kern-\wd0 \kern0,4pt \copy0 \kern-\wd0
     \kern-0,2pt \raise0,3pt \box0 }
\def\pt{\PD{}{t}} 
\newcommand{\aop}{\opr{C}}			
\newcommand{\amat}{\arrfont{C}}			
\newcommand{\arrfont}[1]{\mbox{\sffamily$\textbf{#1}$}}
\newcommand{\cmt}{$//$}				
\newcommand{\comap}[2][]{#1\vartheta_{#2}}	
\newcommand{\conlab}[1]{#1_{\rm{c}}}		
\newcommand{\Czero}{\set{C}^0}			
\newcommand{\diag}{\operatornamewithlimits{diag}}
\newcommand{\dlebnorm}[1]{||#1||_\pode}		
\newcommand{\DP}[2]{{#1\boldsymbol{\cdot}#2}}	
\newcommand{\dss}{\boldsymbol{\mathsf{\Sigma}}}	
\newcommand{\elbox}[2][]{#1{\set{E}}_{#2}}	
\newcommand{\eld}{h}				
\newcommand{\elsiz}[2][]{\eld_{#2}^{#1}}	
\newcommand{\errterm}[1]{\opr{E}_{#1}}		
\newcommand{\enum}{\lavec{e}}                   
\newcommand{\fnum}{\lavec{f}}			
\newcommand{\GL}{\text{\sc{gl}}}		
\newcommand{\Gn}[2][]{#1\xi_{#2}}		
\newcommand{\Gw}[1]{w_{#1}}			
\newcommand{\grad}{\vec{\nabla}}		
\newcommand{\Hmat}{\arrfont{H}}                 
\newcommand{\Hone}{\set{H}^1}			

\newcommand{\ipc}[2]{\left\langle#1,%
#2\right\rangle}				
\newcommand{\isin}[2]{\in\left\{#1,\cdots#2\right\}}
\newcommand{\Jmat}{\boldsymbol\Phi}		
\newcommand{\lavec}[1]{{\boldsymbol#1}}		
\newcommand{\lclose}{\right\lbrack}		
\newcommand{\Leg}[1]{{\rm L}_{#1}}		
\newcommand{\Lop}{\nabla^2}			
\newcommand{\Lmat}{\arrfont{L}}			
\newcommand{\Ltwo}{\set{L}_2}			
\newcommand{\mask}{\boldsymbol{\mathsf{\Pi}}}	
\newcommand{\Mmat}{\arrfont{M}}			
\newcommand{\mult}{\arrfont{W}}			
\newcommand{\mpar}[1]{
}
\newcommand{\mv}[2][]{{\langle{#2}\rangle_{#1}}}
\newcommand{\nopo}[1][\pode]{N_{#1}}		
\newcommand{\oforder}[1]{\mathcal{O}({#1})}	
\newcommand{\opr}[1]{\mathcal{#1}}		
\newcommand{\PD}[2]{\partial_{#2}#1}		
\newcommand{\pdomain}{\set{D}}			
\newcommand{\Pmat}{\arrfont{P}}			
\newcommand{\pn}[2][]{#1x_{#2}}			
\newcommand{\pode}{p}				
\newcommand{\polynomialsset}[1]{\set{V}_{#1}}	
\newcommand{\pwisepolysset}[2]{\polynomialsset{#1,#2}}
\newcommand{\Qmat}{\arrfont{A}}			
\newcommand{\rnum}{\lavec{r}}			
\newcommand{\ropen}{\left\rbrack}		

\newcommand{\sd}{{\scriptscriptstyle\ddots}}
\newcommand{\set}[1]{\mathbb{#1}}		
\newcommand{\setdef}[2]{\left\{#1\left\bracevert#2\right.\right\}}
\newcommand{\si}{\mu}				
\newcommand{\smooth}{\arrfont{S}}	        
\newcommand{\so}{{\scriptscriptstyle1}}
\newcommand{\sz}{{\scriptscriptstyle0}}
\newcommand{\spanop}[2][]{\operatorname{span}^{#1}_{#2}}

\newcommand\tendsto[2]{\xrightarrow[#1\rightarrow#2]{}}
\newcommand{\tr}[1]{\raisebox{0pt}
{$#1$}^\text{\sc{t}}}
\newcommand{\ua}{\vec{\uac}}			
\newcommand{\uac}{c}

\newcommand{\union}[2][]{\bigcup_{#2}^{#1}}
\newcommand{\unum}{\lavec{u}}			
\newcommand{\vnum}{\lavec{v}}			
\newcommand{\uspace}{\set{U}}			
\newcommand{\wnum}{\lavec{w}}			
\newcommand{\uv}[1]{{\vec{\iota}}^{\:#1}}	
\setlength{\marginparsep}{\baselineskip}
\setlength{\marginparwidth}{3.2cm}
\setlength{\parindent}{1em}	

\date{11 May 2005}
\begin{document}
\begin{frontmatter}
%
%
%
\title{Geophysical-astrophysical spectral-element adaptive refinement
(\gaspar): Object-oriented $h$-adaptive code for geophysical fluid dynamics
simulation}
\author[ncar:image]{Duane Rosenberg}
\ead{duaner@ucar.edu}
\author[ncar:image]{Aim\'e Fournier}
\ead{fournier@ucar.edu}
\author[anl]{Paul Fischer}
\ead{fischer@mcs.anl.gov}
\author[ncar:essl]{Annick Pouquet}
\ead{pouquet@ucar.edu}
\address[ncar:image]{Institute for Mathematics Applied to Geosciences}
\address[ncar:essl]{Earth and Sun Systems Laboratory\\ 
National Center for Atmospheric Research\\ 
PO Box 3000, Boulder, Colorado 80307-3000 USA}
\address[anl]{Mathematics and Computer Science Division\\
Argonne National Laboratory, Illinois 60439-4844 USA}
\begin{abstract}
We present an object-oriented geophysical and astrophysical
spectral-element adaptive refinement (\gaspar) code for application to
turbulent flows.
Like most spectral-element codes,
\gaspar\ combines finite-element efficiency with spectral-method accuracy.
It is also designed to be flexible enough for
a range of geophysics and astrophysics applications where turbulence
or other complex multiscale problems arise.
For extensibility and flexibilty the code is designed in an object-oriented
manner.
The computational core is based on spectral-element operators, which are
represented as objects. The formalism accommodates both conforming and
nonconforming elements and their associated data structures for
handling interelement communications in a parallel environment. Many aspects
of this code are a synthesis of existing methods; however, we focus on a 
new formulation of dynamic adaptive refinement (DARe) of nonconforming $h$-type. 
This paper presents the code and its algorithms; we do not
consider parallel efficiency metrics or performance. As a demonstration of the
code we offer several two-dimensional test cases that we propose as
standard test problems for comparable DARe codes. The
suitability of these test problems for 
turbulent flow simulation is considered. 
\end{abstract}
\end{frontmatter}
\setlength{\parskip}{1\parskip}	

\section{Introduction}
\label{sec:intro}

Accurate and efficient simulation of strongly turbulent flows is a prevalent challenge in many atmospheric, oceanic, and
astrophysical applications.
New numerical methods are needed to investigate such
flows in the parameter regimes that interest the geophysics communities.
Turbulence is one of the
last unsolved classical physics problems, and such flows
today form the focus of numerous investigations.
They
are linked to many issues in the geosciences, for example , in meteorology,
oceanography, climatology, ecology, solar-terrestrial interactions,
and solar fusion, as well as dynamo effects, specifically, magnetic-field generation in
cosmic bodies by turbulent motions.
Nonlinearities prevail in turbulent flows when the
Reynolds number $\Rn$
is large.
The number of degrees of freedom (\dof) increases as
$\Rn^{9/4}$ for $\Rn\gg1$ in the Kolmogorov (1941) framework.
\mpar{cite specific K41 and $\Rn^{9/4}$ references.}
For aeronautic flows often $\Rn>10^6$, but for geophysical flows often $\Rn\gg10^8$;
\cite{MK00, ES04}
for this and other reasons the ability to probe large $\Rn$, and 
to examine in detail the
large-scale behavior of turbulent flows depends
critically on the numerical ability to resolve a large number of
spatial and temporal scales.

\mpar{One intriguing observation concerning turbulent flows resides in the
departure from normality in the probability distribution
functions. The origin of these fat wings is not understood;
one does not know yet what structures are key to our
understanding turbulent flows' statistical properties, {\it e.g.,}
vortex sheets, spirals or filaments,
shocks or fronts, blobs, plumes or tetrads, knots, helices, tubes or
arches.
The link between structures within the flow and the observed
non-Gaussian statistics of the velocity gradient matrix is the basis
for the notion of intermittency which plays a role {\it e.g.,} in
reactive flows, in convective plumes, in localized propagating fronts,
in combustion or in solar-corona heating.}

Theory demands that computations of nonconvective turbulent flows
reflect a clear separation between the energy-containing 
and dissipative scale ranges.
\mpar{Cite what theory demands this.}
Uniform-grid convergence
studies on 3D compressible-flow simulations show that in order to achieve the desired
scale separation, 
uniform grids must contain at least $2048^3$ cells \cite{sytine2000}.
Today such computations can barely be accomplished.
A pseudo-spectral Navier-Stokes code on a grid of $4096^3$ uniformly spaced
points has been run on the Earth Simulator \cite{isihara},
but the Taylor Reynolds number ($\propto\sqrt{\Rn}$) is still no more than 
$\approx700$, very far from what is required for most geophysical flows.
Similar scale 
separation is required in computations of convective-turbulent flows.
However, if the flow's significant structures are indeed sparse,
so that their dynamics can be followed accurately even if they are
embedded in random noise, then dynamic adaptivity seems to offer a
means for achieving otherwise unattainable large $\sqrt{\Rn}$ values.

We have developed a dynamic geophysical and astrophysical spectral-element adaptive refinement
(\gaspar) code for simulating and studying turbulent phenomena.
Several properties of spectral-element methods (SEMs)
make them desirable for direct numerical simulation
(DNS) or large-eddy simulation (LES) of geophysical turbulence. 
Perhaps most significant is the 
fact that SEM are inherently minimally diffusive and dispersive. This
property is clearly important when trying to simulate high $\Rn$ number (low
viscosity) flows that characterize turbulent behavior. Also, because
SEMs use a finite element as a fundamental description of the geometry 
\cite{patera84}, they can be used in very efficient high-resolution turbulence studies in 
domains with complicated boundaries. Such discretizations also enable SEMs
to be naturally parallelizable (\cite[\eg,][]{FTT2004}), an important feature when simulating
flows at high $\Rn$ with many \dof\ involving multiple
spatial and temporal scales. Equally important, spectral element methods not only enjoy
spectral convergence properties when the solution is smooth but are also effective
when the solution is not smooth \cite{DFM2002} (see \eqref{e:errorbound} below).

\begin{figure}
\begin{center}
\includegraphics[scale=.60]{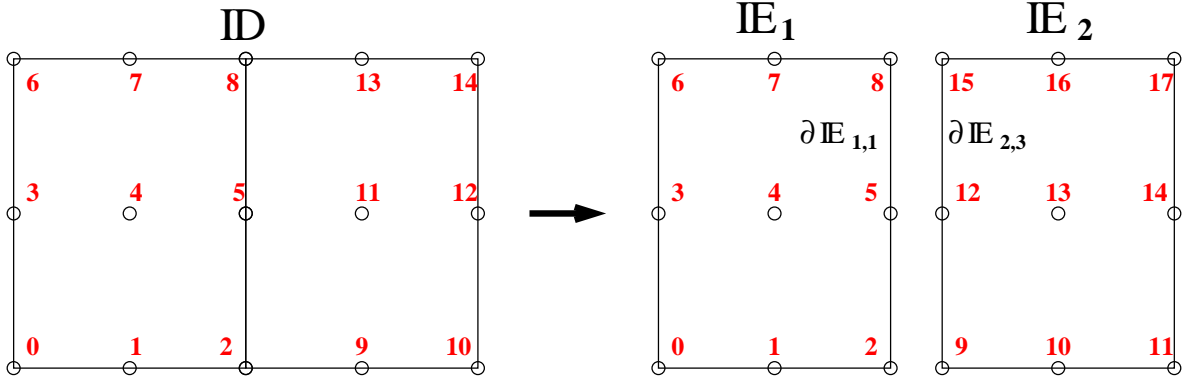}
\caption{Schematic of a conforming degree $\pode=2$ mesh showing the mapping of
global (i.e., unique) \dof\ in the problem domain $\pdomain$ (left) to local (i.e., redundant) \dof\ in the elements $\elbox{k}$ (right).
Edge subscripts give element index $k$ and edge index from $s=0$ (south
edge) counterclockwise to $s=3$. 
The element $\elbox{1}$ is bounded at the east by the edge $\partial\elbox{1,1}$ and
$\elbox{2}$ is bounded at the west by the edge $\partial\elbox{2,3}=\partial\elbox{1,1}$.
The interface matching condition occurs by simple assignment leading to a Boolean assembly matrix $\conlab{\Qmat}$.
\vspace{.2\baselineskip}	
}
\label{fig:conf_config}
\end{center}
\end{figure}

In the case of SEMs, conforming adaptive methods
(where entire element edges geometrically coincide, as in \fig{fig:conf_config})
are gradually being replaced by nonconforming adaptive methods. One reason is that
mesh generation for conforming methods is complicated when attempting to resolve
local flow features \cite{ronquist96}. Another reason is that adaptive 
conforming meshes can lead to high-aspect-ratio elements that can cause difficulties
for a linear solver \cite{fischer_kruse_2002}. Moreover, the fact that nonconforming elements
can better localize mesh refinement implies that the computational cost among all
elements can be reduced \cite{krusediss,orszag80}; with conforming 
adaption, however, local refinement regions may extend out to where refinement is not
dictated by local features of interest within the solution.
\mpar{Why is \cite{orszag80} cited here?}

\begin{figure}
\begin{center}
\includegraphics[scale=.60]{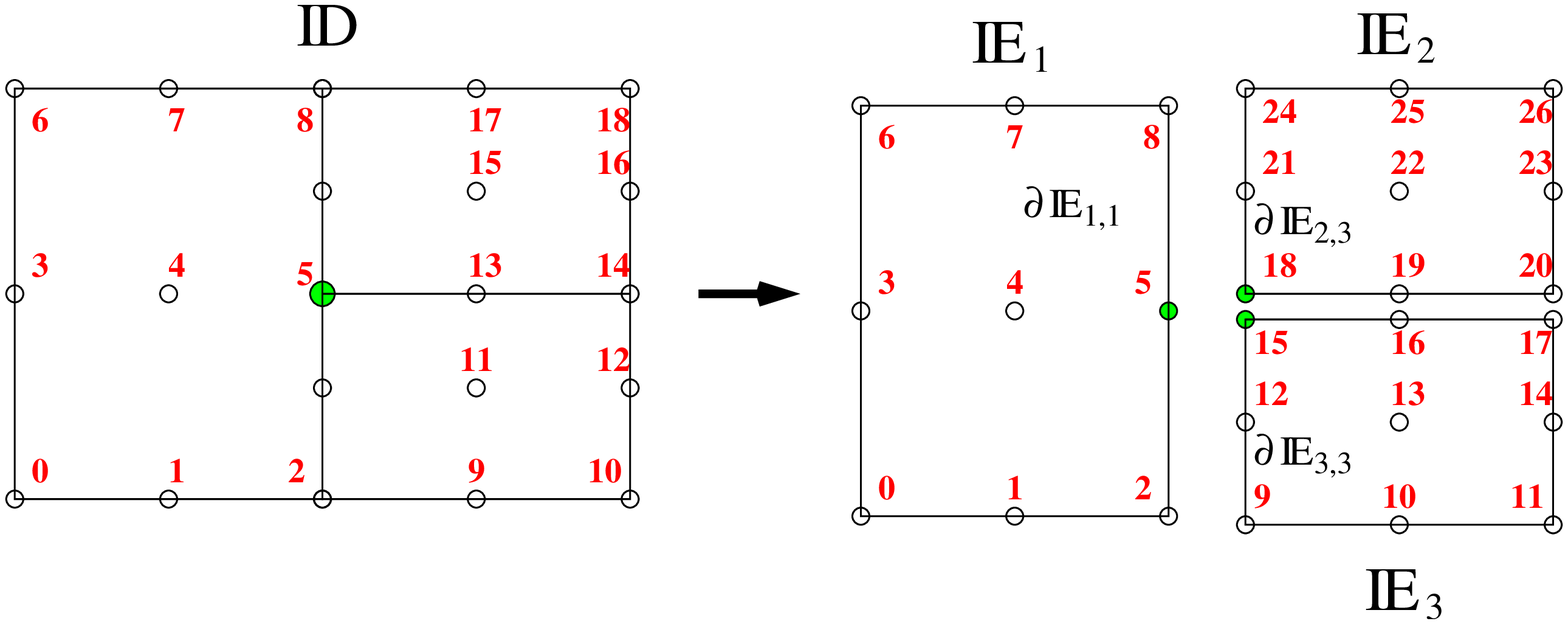}
\end{center}
\caption{As in \fig{fig:conf_config} but for a geometrically \emph{non}conforming (but functionally conforming---each element has $\pode=2$) mesh.
Here $\elbox{2}$ and $\elbox{3}$ are bounded at the west by ``child''
edges $\partial\elbox{2,3}$ and $\partial\elbox{3,3}$, and
$\elbox{1}$ is bounded at the east by the ``parent'' edge $\partial\elbox{1,1}=\partial\elbox{2,3}\bigcup\partial\elbox{3,3}$.
The mapping occurs by means of an interpolation of global \dof\ from the function space associated with the parent edge
onto the union of those associated with the child edges, which contains the function space of the parent.
}
\label{fig:nonconf_config}
\end{figure}

Nonconforming element discretizations can be \emph{geometrically} and/or \emph{functionally}
nonconforming. In the former case (\fig{fig:nonconf_config}), neighboring elements have
boundaries that do not
coincide; in the latter, the polynomial expansion degree $\pode$
in neighboring elements can differ.
Several SEM
researchers have adopted the discretization method that simultaneously
alters element size $h$ and configuration ($h$-refinement) {\em and} the polynomial
degree $\pode$ in neighboring elements ($p$-refinement), providing for a so-called
$h$-$p$-refinement strategy.
The \emph{mortar element method} (MEM)
\cite{anagnostou89,bernardi89,maday89} is a nonconforming discretization method that
uses a variational formulation to minimize the Lebesgue $\Ltwo$ norm of the discontinuous jump
across nonconforming spectral-element boundaries. In
a number of recent applications, this method has been used in unstructured
(static, nonadaptive) simulations of turbulence \cite{henderson95} and for
ocean dynamics \cite{iskandarani2003,levin2000}.
This method has been shown to produce optimal
convergence in the incompressible Stokes equations solution \cite{belgacem2000},
and it has been demonstrated experimentally to produce excellent results when used
as a basis for adaptive mesh refinement \cite{mavriplis94}.

\mpar{This might be a better order of the same sentences in this \P.}
Nonconforming $h$-$p$ adaptive methods using MEM have been developed
for studying turbulent phenomena \cite{henderson99}, ocean modeling
\cite{levin2000}, flame front deformation \cite{hfeng2002},
electromagnetic scattering \cite{KoWoHu02}, wave propagation
\cite{CGFMQ02}, seismology \cite{ChCaVi03} and other topics.
However, MEM for $p$-type refinement has been cited as
causing instabilities in flow calculations \cite{ronquist96};
for this reason the interpolation-based method was developed.
Also, in most flows of interest to us, it is the different
scales' interaction that determines not only the structures
that form but also their statistics and time evolution.
This suggests that reasonably high-order approximations
are required in each element during much of the evolution.
Thus, in the present work, using the fact that the
numerical solution order is relatively high, we restrict
ourselves to a nonconforming $h$-refinement strategy only and
use an {\it interpolation-based} scheme to maintain continuity
between nonconforming elements \cite{fischer_kruse_2002,krusediss}.
(Throughout the remainder of this paper ``nonconforming''
will refer to geometrically nonconforming elements,
keeping the expansion degree $\pode$ fixed in all elements.)
Researchers have recently used this nonconforming treatment as a basis for performing
DARe \cite{StCyr2004b};
however, to the best of our knowledge, our implementation of this
interpolation scheme in the fully dynamic adaptivity context is unique.

Our purpose in this paper is to describe \gaspar and, in particular, 
the procedures used in the DARe technique. We first
describe (\sr{sec:discretization}) SEM discretization on a particular class of problems 
and introduce many of the required formulas, operators, and so forth.
We discuss (\sr{sec:connectivity}) the type of nonconforming discretization allowed 
and the way in which continuity is maintained between elements sharing nonconforming interfaces.
We explain the linear-solver details in \sr{sec:linear_solver}.
In 
\sr{sec:adaptivity} we consider how nonconforming edges are identified 
and how neighboring elements are found, we and present the rules for refinement and coarsening.
In \sr{sec:data_comm} we present the gather/scatter
matrix that facilitates communication of edge data. A priori error estimators
are discussed in \sr{sec:error} as a means to provide DARe criteria. 
Then, in \sr{sec:test_problems} we consider two test-problem classes: 
solutions of the heat equation and linear advection-diffusion
equation (\sr{sec:linadvecttest}), which demonstrate feature tracking of smooth and
isolated features governed by dynamics of linear systems; and the Burgers
equation, that tests the ability of DARe to capture and track well-defined 
and reasonably sharp structures (\sr{sec:burgerstest}) that arise from nonlinear dynamics.
In \sr{sec:conclusion} we offer some conclusions on our current work, 
as well as some comments on the application of
\gaspar\ to geophysical turbulence simulations.

\section{Methodology}

Because we wish to focus on the DARe methodology of \gaspar,
we concentrate on the simplest multidimensional nonlinear
equation that encompasses many of the difficulties in simulating
fluid turbulence. Thus in this section we present the discrete form for 
the 2D Burgers equation and its variants.
In turn we present the spatial discretization using a variational
formulation and the spectral-element operators that arise. We then
present the time discretizations.

\subsection{Discretization of the dynamics}
\label{sec:discretization}

The equations considered in this work derive from the
$d$-dimensional advection-diffusion equation for velocity
$\u(\x,t)$:
\be
\pt \u +\DP{\ua}{\grad\u} = \nu \Lop\u,
\label{eq:burgers}
\ee
where $\nu\propto\Rn^{-1}$ is the kinematic viscosity.
This is to be solved in a spatiotemporal domain $(\x,t)\in\pdomain\times\left]0,t_{\rm f}\right]$
subject to the boundary and initial conditions
\begin{align}
\u(\x,t)&=\vec{b}(\x,t)\quad\text{for}\quad(\x,t)\in\partial\pdomain\times\left]0,t_{\rm f}\right],
\label{e:BoCo}\\
\qquad\u(\x,0)&=\u_{\rm{i}}(\x)\quad\text{for}\quad\x\in\pdomain.
\label{e:InCo}
\end{align}
In this work, $\ua$ may be $\u$, so that \eq{eq:burgers} is the Burgers equation, or
$\ua=\vec{c}(t)$, a prescribed uniform linear advection velocity.
\mpar{We want to specify the spaces for u and the test function, v, to set up the
fact that $\Czero$ is all that is required for maintaining them. We can then show, if
we want, that the error with P is exponential (i.e., spectral).
 When we establish that continuity is all that is required, we can make a statement
 about the fact that conforming elements will take care of this, then we get to
 nonconforming treatment (which is really a conforming treatment in that 
 the conforming Lag. bases are evaluated on the child-element interfaces...}

\subsubsection{Variational approach to spatial discretization}
\label{sec:variational_form}

Define the function space
$$
\uspace_{\vec{b}}\equiv%
\setdef{\u=\sum_{\si=1}^du^\si\uv{\si}\,}{u^\si\in\Hone(\pdomain)\thickspace\forall\mu\quad\&\quad\u=\vec{b}\thickspace%
\text{on}\thickspace\partial\pdomain},
$$
where $\uv{\si}$ denotes the Cartesian unit vectors and
$$
\Hone(\pdomain)\equiv\setdef{f\,}{f\in\Ltwo(\pdomain)\quad\&\quad\PD{f}{x^\si}\in\Ltwo(\pdomain)\thickspace\forall\mu}.
$$
Then the discretization of \eq{eq:burgers} is based on the following variational form:
Find the trial function $\u(\cdot,t)\in\uspace_{\vec{b}}$ such that for any test function $\v\in\uspace_{\vec{0}}$
\be
\ipc{\v}{\pt\u}+\ipc{\v}{\aop\u}=-\nu\ipc{\tr{\grad\v}}{\grad\u}, 
\label{eq:weak}
\ee
where $\aop\equiv\DP{\ua}{\grad}$ is the advection operator and the inner product is
\eqref{e:innpro}.
Note that the treatment of \eqref{e:InCo} will not be made explicit but may be easily inferred from
the general discussion.

We assume that $\pdomain$ can be partitioned as in \eqref{e:domdecom}.
(See the appendix for the complete mathematical details.)
To discretize \eq{eq:weak}, we adopt a Gauss-Lobatto-Legendre
(GLL) basis \eq{e:limappedv}, that is, expand $u^\si$ and $v^\si$
using \eqref{projmapped}.
Inserting these expansions into \eq{eq:weak}, we arrive at the semi-discrete ODE system
problem: Find the numerical solution
$\u_{\rm{n}}(\cdot,t)=\tr{\smash{\vec{\lavec{\phi}}}}\unum(t)\in\opr{P}_{\lavec{\eld},\vec\pode}\uspace_{\vec{b}}$
such that for all
$\v=\tr{\smash{\vec{\lavec{\phi}}}}\vnum\in\opr{P}_{\lavec{\eld},\vec\pode}\uspace_{\vec{0}}$,
\be
\tr{\vnum} \Mmat\frac{\d\unum}{\d t}+\tr{\vnum}\amat\unum =-\nu\tr{\vnum}\Lmat\unum
\label{eq:vdiscrete_burgers}
\ee
collocated
at $K\nopo^d$ mapped Lagrange node points \eqref{e:npmapped}%
, where 
$\Mmat=\diag_k\Mmat_k$, $\amat=\diag_k\amat_k$, and $\Lmat=\diag_k\Lmat_k$
are the unassembled block-diagonal mass matrix, linear or nonlinear advection matrix
\cite[\cf][Ch.\ 6]{DFM2002}, and diffusion matrix, respectively.
The respective $d\nopo^d\times d\nopo^d$ matrix blocks locally indexed to element $\elbox{k}$ are
\be
M_{\vec{\jmath},\vec{\jmath}\,';k}^{\si,\si'}\equiv%
\ipc{\vec{\phi}_{\vec{\jmath},k}^{\si}}{\vec{\phi}_{\vec{\jmath}\,',k}^{\si'}}_\GL=
\delta_{\vec{\jmath},\vec{\jmath}\,'}\delta^{\si,\si'}\Gw{\vec{\jmath},k},
\label{e:massop}
\ee
including the quadrature weights \eqref{e:GLLqw},
\be
C_{\vec{\jmath},\vec{\jmath}\,';k}^{\si,\si'}\equiv%
\ipc{\vec{\phi}_{\vec{\jmath},k}^{\si}}{\aop\vec{\phi}_{\vec{\jmath}\,',k}^{\si'}}_\GL=
\delta^{\si,\si'}\Gw{\vec{\jmath},k}\DP{\ua_{\vec{\jmath},k}}{\grad}%
\phi_{\vec{\jmath}\,',k}(\pn[\vec]{\vec{\jmath},k}),
\label{e:advop}
\ee
where $\ua_{\vec{\jmath},k}(t)\equiv\ua(\pn[\vec]{\vec{\jmath},k},t)$, and
$$
L_{\vec{\jmath},\vec{\jmath}\,';k}^{\si,\si'}\equiv\ipc%
{\grad\vec{\phi}_{\vec{\jmath},k}^{\si}}%
{\grad\vec{\phi}_{\vec{\jmath}\,',k}^{\si'}}_\GL=
\delta^{\si,\si'}\sum_{\vec{\jmath}\,''\in\set{J}}%
\Gw{\vec{\jmath}\,'',k}\DP%
{\grad\phi_{\vec{\jmath},k}(\pn[\vec]{\vec{\jmath}\,'',k})}%
{\grad\phi_{\vec{\jmath}\,',k}(\pn[\vec]{\vec{\jmath}\,'',k})}.
$$
To compute $\grad\phi_{\vec{\jmath},k}$, one differentiates \eqref{e:lagrinteex} and uses \eqref{e:limapped}.
The weak diffusion matrix for deformed quadrilaterals or deformed cubes
(nonlinear maps $\comap[\vec]{k}$) can also be constructed
(\eg, \cite{DFM2002}). While these elements are supported in \gaspar,
they are not necessary for the present discussion.

Restricting ourselves for the moment to the $k$th
element $\elbox{k}$, we see that \eq{eq:vdiscrete_burgers} must hold for
the restriction $\v|_{\elbox[\bar]{k}}=\tr{\smash{\vec{\lavec{\phi}}}}_k\vnum_k$ of $\v$ to any $\elbox{k}$, so that
the unassembled ODE for $\u_{\rm{n}}|_{\elbox[\bar]{k}}=\tr{\smash{\vec{\lavec{\phi}}}}_k\unum_k$ is
\be
\Mmat_k\frac{\d\unum_k}{\d t}+\amat_k\unum_k=-\nu\Lmat_k\unum_k.
\label{eq:discrete_burgers}
\ee
Strictly speaking, \eqref{eq:discrete_burgers} is true only after being assembled in the
manner discussed in \sr{sec:global_assembly}.
Continuity of $\u_{\rm{n}}$ across all elements is a sufficient condition for maintaining $u_{\rm{n}}^\si\in\Hone(\pdomain)$. We allow
two element configuration types, conforming and nonconforming, as illustrated in 
Figs.\ \ref{fig:conf_config} and \ref{fig:nonconf_config}, respectively.
Conforming discretizations enforce continuity simply by assigning the same
$\u_{\rm{n}}$ values to the coinciding node points $\pn[\vec]{\vec{\jmath},k}=\pn[\vec]{\vec{\jmath}\,',k'}$
along element edges $\partial\elbox{k,s}=\partial\elbox{k',s'}$.
\mpar{DLR: Bring up the gather-scatter method here?}
In the nonconforming case $\partial\elbox{k,s}\subsetneq\partial\elbox{k',s'}$
and functional continuity cannot be accomplished
by simple assignment because most node points are not coinciding; thus, we use
an interpolation-based scheme to enforce continuity along a nonconforming interface. 
This scheme is the subject of \sr{sec:connectivity}.

\subsubsection{Time discretization}
\label{sec:time_discrete}
While there are many time discretization schemes (\eg,\cite{canuto}),
we restrict ourselves to semi-implicit multistep methods.
In all these methods the diffusion is solved fully implicitly while the
time-derivative is approximated using a backward-difference formula
(BDF) of order $M_{\rm bdf}$ (see \cite{DFM2002,KS99}) and the
advection term is approximated by using an explicit 
extrapolation-based method (Ext) of order $M_{\rm{ext}}$ \cite{KIO91}.
Then, to integrate \eq{eq:discrete_burgers} from time $t^{n-1}$ to time
$t^n$, one has
\be
\Hmat_k^n\unum_k^n= 
\sum_{m=n-M_{\rm bdf}}^{n-1}\beta_{\rm bdf}^{m,n}\Mmat_k^m\unum_k^m
-\sum_{m=n-M_{\rm{ext}}}^{n-1}\beta_{\rm{ext}}^{m,n}\amat_k^m\unum_k^m,
\label{eq:time_desc_burgers}
\ee
where
\be
\Hmat_k^n\equiv\beta_{\rm bdf}^{n,n}\Mmat_k^n+\nu\Lmat_k^n
\label{eq:helmholtzop}
\ee
is a discrete spectral-element 
Helmholtz operator. 
Although the matrices $\Lmat_k$ and $\Mmat_k$ in \eq{eq:discrete_burgers} were $t$-independent,
they are time-indexed in \eqref{eq:time_desc_burgers} because
DARe will, in general, reconfigure the partition \eqref{e:domdecom} over time.
For this reason the coefficients $\beta^{m,n}$ are computed for each $t^n$
as in the traditional schemes cited except that the timestep $\varDelta{t}$ may vary
with $m$ as the smallest spectral-element diameter $\elsiz{}\equiv\min_k\elsiz{k}$ \eqref{e:elsiz} changes.
As discussed below, $\u_{\rm{n}}^{\,n}$ continuity is maintained during
the linear solve of \eqref{eq:time_desc_burgers} assembled over $k$
for the solution $\unum^n$.
Because the matrix $\Hmat^n\equiv\diag_k\Hmat_k^n$ is symmetric positive-definite (SPD),
provided that $\u_{\rm{n}}^{\,n}$ is restricted to $\uspace_{\vec{0}}$,
the solution of the assembled \eq{eq:time_desc_burgers} is obtained by using a preconditioned 
conjugate gradient
(PCG, see \cite{shewchuck94,WeisLQ2}) algorithm to invert $\Hmat^n$ at the time step to $t^n$. 
In the presence of nonconforming elements, care must be exercised to ensure that the
search directions in the PCG algorithm are in $\uspace_{\vec{0}}$. We consider 
appropriate modifications to PCG in \sr{sec:linear_solver}.

\subsection{Implications for code design}

The fully discretized advection-diffusion equation \eq{eq:time_desc_burgers} suggests
several issues that we have taken into consideration when designing the
code. First, all geometric (mesh) information is separated from all other
code objects, since element type information can be encoded
easily into the objects that require this distinction. Second, the solution
data must be available at multiple time levels, so it is reasonable to 
provide this information in a data structure. Thus we are led to 
the notion of {\it element} and {\it field} objects. 
The former contains all geometric information in $d$ dimensions, including the values and 
tensor-product ordering of the Gauss-quadrature nodes \eqref{e:npmapped} and weights \eqref{e:GLLqw}.
The element object also contains neighbor-list information 
and an encoding of the hierarchical element refinement level
$\propto\log_{\half}\elsiz{k}$ of each element $\elbox{k}$. The field object contains 
the data $\unum^m$ representing the physical quantity of interest at all required time levels $t^m$.
 
We note that, while we concern ourselves with the advection-diffusion equation \eqref{eq:burgers} here, the code
will allow any equation that is discretized by using a Lagrangian basis on up to two meshes
(important, for example, when discretizing the Navier-Stokes equations by using a
formulation coupling the polynomial spaces $\polynomialsset{\pode}$-$\polynomialsset{\pode-2}$ \cite{maday92}).
The basis
functions and the 1D derivative matrices and Gauss-quadrature nodes
\eqref{e:npmapped1} and weights \eqref{e:GLLqwmapped} are encapsulated within
basis classes (objects), and the SEM operators such as (\ref{e:massop},\ref{e:advop},\ref{eq:helmholtzop}) above are
constructed as objects that contain pointers to the basis objects and to a local element
object. Generally the $d$-dimensional SEM
operators are not stored but are
constructed from a tensor product \emph{application}
of the relevant 1D objects. High-level objects 
encapsulate the solution of \eq{eq:discrete_burgers} or other equations, such as the Navier-Stokes
equations, and have
common interfaces that allow the equations to take a single time integration step; in 
other words, all
high-level equation solver classes are used in the same way. 
All the higher-level objects are constructed with ``smart'' arrays or linked lists 
of elements and fields that are independent of the objects that solve the equations. 
Hence, the classes that handle DARe and enforce continuity between elements 
being solved. The high-level objects also contain linked lists of SEM operators that \emph{do}
depend on the equation being solved. 

\subsection{Local application of spectral-element operators}
\label{sec:multiple_subdomains}

In general, global meshes consist of multiple-element grids.
In this subsection, we make explicit the form of the SEM operators 
implemented in the code that ensure the solution continuity 
across element interfaces.

\subsubsection{Continuity between nonconforming elements}
\label{sec:connectivity} 

In a conforming treatment, $\Czero$ continuity is maintained between elements
by ensuring that the function values on the coinciding nodes are the same.
The matching condition, then, consists of expressing the $N_{\rm{g}}$ global (unique)
\dof\ $\unum_{\rm{g}}$ in terms of the local (redundant) \dof\ as $d\nopo^d$-vectors $\unum_k$,
$k\isin{1}{K}$.
Generally $N_{\rm{g}}<Kd\nopo^d$.
This mapping can be accomplished by using a $Kd\nopo^d\times N_{\rm{g}}$
Boolean assembly matrix $\conlab{\Qmat}$, such that
\begin{equation}
\unum = \conlab{\Qmat}\unum_{\rm{g}}.
\label{eq:u1Qug}
\end{equation}

For example, if our mesh partition is that in 
\fig{fig:conf_config}, we can easily see that \eqref{eq:u1Qug} takes the following explicit
form (suppressing zero-valued and $\si>1$ blocks):
{
\renewcommand{\arraystretch}{.3}
\setlength{\arrayrulewidth}{.1pt}
$$
\setcounter{MaxMatrixCols}{15}
\unum=\begin{pmatrix}
u_0\\\vdots\\u_{17}\end{pmatrix}=
\begin{pmatrix}
u_{0,1}\\\vdots\\ u_{8,1}\\\hline
u_{0,2}\\\vdots\\ u_{8,2}\end{pmatrix}=\begin{pmatrix}
   \so&  \sz&  \sz&  &&&&&&&&&&&\\
   \sz&  \so&  \sz&  &&&&&&&&&&&\\
   \sz&  \sz&  \so&  &&&&&&&&&&&\\
   &&&\so&  \sz&  \sz&  &&&&&&&&\\
   &&&\sz&  \so&  \sz&  &&&&&&&&\\
   &&&\sz&  \sz&  \so&  &&&&&&&&\\
   &&&&&&\so&  \sz&  \sz&  &&&&&\\
   &&&&&&\sz&  \so&  \sz&  &&&&&\\
   &&&&&&\sz&  \sz&  \so&  &&&&&\\\hline
   \sz&  \sz&  \so&  &&&&&&\sz&  \sz&  &&&\\
   \sz&  \sz&  \sz&  &&&&&&\so&  \sz&  &&&\\
   \sz&  \sz&  \sz&  &&&&&&\sz&  \so&  &&&\\
   &&&\sz&  \sz&  \so&  &&&&&\sz&  \sz&  &\\
   &&&\sz&  \sz&  \sz&  &&&&&\so&  \sz&  &\\
   &&&\sz&  \sz&  \sz&  &&&&&\sz&  \so&  &\\
   &&&&&&\sz&  \sz&  \so&  &&&&\sz&  \sz\\
   &&&&&&\sz&  \sz&  \sz&  &&&&\so&  \sz\\
   &&&&&&\sz&  \sz&  \sz&  &&&&\sz&  \so\\
\end{pmatrix}\begin{pmatrix}
u_{{\rm{g}},0}\\\vdots\\u_{{\rm{g}},14}\end{pmatrix}.
\setcounter{MaxMatrixCols}{10}
$$
}

In practice, $\conlab{\Qmat}$ is never formed explicitly but is instead
{\it applied}. The matrix $\conlab{\Qmat}$, which assigns global \dofs\
to the local \dofs, is also called a {\it scatter} matrix.
Its transpose $\conlab{\tr{\Qmat}}$ performs the associated {\it gather} operation.

In our nonconforming treatment,
we use an interpolation-based method for establishing the interface matching 
conditions between neighboring (nonconforming) elements. The underlying function space
$\uspace_{\vec{b}}$ (\sr{sec:discretization}) remains the same, unlike
in MEM. Consider the nonconforming mesh configuration depicted in
\fig{fig:nonconf_config}.
%
For the moment denote the global 
nodes, those nodes residing on the east {\it parent} edge $\partial\elbox{1,1}$,
by $\pn[\vec]{{\rm{g}},i}$, $i\in\{2,5,8\}$,
and denote the nodes on the west {\it child} edges, $\partial\pdomain_{2,3}$ and
$\partial\pdomain_{3,3}$ by $\vec{x}_j$, $j\in\{9,12,15,18,21,24\}$.
We explain elsewhere how in the weak formulation of \eqref{eq:burgers} or other PDEs, the
test function factors $\phi_{{\rm{g}},i}(\vec{x})$ arise that are continuous across the global domain $\pdomain$,
and interpolate from the global node values $\pn[\vec]{{\rm{g}},i}$
\cite{FRP2005}.
A continuous solution is then found in
$\spanop{i}\phi_{{\rm{g}},i}$ by projecting there from the space of trial 
functions $\phi_j(\vec{x})$ that interpolate from the local node values $\vec{x}_j$
but that do not need to be globally continuous.
The matrix block $\phi_{{\rm{g}},i}(\vec{x}_j)$ of $\Qmat$ thus
generalizes the Boolean scatter matrix used in the conforming-element formulation and
accommodates both conforming {\it and} nonconforming elements.
It is convenient to factor $\Qmat
=\Jmat\conlab{\Qmat}
$, where $\Jmat$ is the interpolation matrix from global to local \dofs\
and $\conlab{\Qmat}$ is locally conforming.

For example, the explicit form of \eqref{eq:u1Qug} for the nonconforming assembly matrix
$\Qmat=\Jmat\conlab{\Qmat}$ for the mesh
illustrated in \fig{fig:nonconf_config} is (suppressing zero-valued and $\si>1$ blocks) 
{
\renewcommand{\arraystretch}{.3}
\setlength{\arrayrulewidth}{.1pt}
\begin{align*}
\unum&=\begin{pmatrix}u_0\\\vdots\\u_{26}\end{pmatrix}=
\begin{pmatrix}u_{0,1}\\\vdots\\ u_{8,1}\\\hline
u_{0,2}\\\vdots\\ u_{8,2}\\\hline u_{0,3}\\\vdots\\ u_{8,3}\end{pmatrix}=\left(
\begin{array}{rrlrrlrrlrrrrrrrrrr}     
   \so&  \sz&  \sz&  &&&&&&&&&&&&\\
   \sz&  \so&  \sz&  &&&&&&&&&&&&\\
   \sz&  \sz&  \so&  &&&&&&&&&&&&\\
   &&&\so&  \sz&  \sz&  &&&&&&&&&\\
   &&&\sz&  \so&  \sz&  &&&&&&&&&\\
   &&&\sz&  \sz&  \so&  &&&&&&&&&\\
   &&&&&&\so&  \sz&  \sz&  &&&\\
   &&&&&&\sz&  \so&  \sz&  &&&\\
   &&&&&&\sz&  \sz&  \so&  &&&\\\hline
   \sz&  \sz&  \so&  &&&&&&\sz&  \sz&  &&&&&&&\\
   \sz&  \sz&  \sz&  &&&&&&\so&  \sz&  &&&&&&&\\
   \sz&  \sz&  \sz&  &&&&&&\sz&  \so&  &&&&&&&\\
   \sz&  \sz&\scriptscriptstyle3/8& \sz&  \sz&\scriptscriptstyle3/4& \sz&  \sz&\scriptscriptstyle-1/8&&&\sz&  \sz&  &&&&&\\
   \sz&  \sz&  \sz&  \sz&  \sz&  \sz&  \sz&  \sz&    \sz&&&\so&  \sz&  &&&&&\\
   \sz&  \sz&  \sz&  \sz&  \sz&  \sz&  \sz&  \sz&    \sz&&&\sz&  \so&  &&&&&\\
   &&&\sz&  \sz&  \so&  &&&&&&&\sz&  \sz&  &&&\\
   &&&\sz&  \sz&  \sz&  &&&&&&&\so&  \sz&  &&&\\
   &&&\sz&  \sz&  \sz&  &&&&&&&\sz&  \so&  &&&\\\hline
   &&&\sz&  \sz&  \so&  &&&&&&&\sz&  \sz&  &&&\\
   &&&\sz&  \sz&  \sz&  &&&&&&&\so&  \sz&  &&&\\
   &&&\sz&  \sz&  \sz&  &&&&&&&\sz&  \so&  &&&\\
   \sz&  \sz&\scriptscriptstyle-1/8&\sz&  \sz&\scriptscriptstyle3/4& \sz&  \sz&\scriptscriptstyle3/8& &&&&&&\sz&  \sz& &\\
   \sz&  \sz&  \sz&  \sz&  \sz&  \sz&  \sz&  \sz&  \sz&  &&&&&&\so&  \sz& &\\
   \sz&  \sz&  \sz&  \sz&  \sz&  \sz&  \sz&  \sz&  \sz&  &&&&&&\sz&  \so& &\\
   &&&&&&\sz&\sz&\so&&&&&&&&&\sz&  \sz\\
   &&&&&&\sz&\sz&\sz&&&&&&&&&\so&  \sz\\
   &&&&&&\sz&\sz&\sz&&&&&&&&&\sz&  \so
\end{array}\right)\begin{pmatrix}
u_{{\rm{g}},0}\\\vdots\\ u_{{\rm{g}},18}\end{pmatrix}\\
&=\left(\begin{array}{rrlrrlrrlrrrrrrrrr}
\so&&&&&&&&&&&&&&&&\\
&&&&&\sd&&&&&&&&&&&\\
&&&&&&&&&&\so&&&&&&\\
\sz&\sz&\scriptscriptstyle3/8&\sz&\sz&\scriptscriptstyle3/4&\sz&\sz&\scriptscriptstyle-1/8&&&&&&&&\\
\sz&\sz&\sz&\sz&\sz&\sz&\sz&\sz&\sz&&&\so&&&&&\\
\sz&\sz&\sz&\sz&\sz&\sz&\sz&\sz&\sz&&&&\sd&&&&\\
&&&&&&&&&&&&&\so&&&\\
\sz&\sz&\scriptscriptstyle-1/8&\sz&\sz&\scriptscriptstyle3/4&\sz&\sz&\scriptscriptstyle3/8&&&&&&&&\\
\sz&\sz&\sz&\sz&\sz&\sz&\sz&\sz&\sz&&&&&&\so&&\\
\sz&\sz&\sz&\sz&\sz&\sz&\sz&\sz&\sz&&&&&&&\sd&\\
&&&&&&&&&&&&&&&&\so
\end{array}\right)\left(\begin{array}{rrrrcrrrrrrrrrrrrrr}
\so&&&&&&&&&&&&&&&&&&\\
&&&&\sd&&&&&&&&&&&&&&\\
&&&&&&&&\so&&&&&&&&&&\\\hline
\sz&\sz&\so&&&&&&&\sz&\sz&&&&&&&&\\
\sz&\sz&\sz&&&&&&&\so&\sz&&&&&&&&\\
\sz&\sz&\sz&&&&&&&\sz&\so&&&&&&&&\\
&&&&&&&&&&&\so&&&&&&&\\
&&&&&&&&&&&&\so&&&&&&\\
&&&\sz&\sz&\so&&&&&&&&\sz&\sz&&&&\\
&&&\sz&\sz&\sz&&&&&&&&\so&\sz&&&&\\
&&&\sz&\sz&\sz&&&&&&&&\sz&\so&&&&\\\hline
&&&\sz&\sz&\so&&&&&&&&\sz&\sz&&&&\\
&&&\sz&\sz&\sz&&&&&&&&\so&\sz&&&&\\
&&&\sz&\sz&\sz&&&&&&&&\sz&\so&&&&\\
&&&&&&&&&&&&&&&\so&&&\\
&&&&&&&&&&&&&&&&\so&&\\
&&&&&&\sz&\sz&\so&&&&&&&&&\sz&\sz\\
&&&&&&\sz&\sz&\sz&&&&&&&&&\so&\sz\\
&&&&&&\sz&\sz&\sz&&&&&&&&&\sz&\so
\end{array}\right)\begin{pmatrix} 
u_{{\rm{g}},0}\\\vdots\\ u_{{\rm{g}},18}\end{pmatrix}.
\end{align*}
}

Note that the $\Qmat$ entries corresponding to the child-node rows (see \fig{fig:nonconf_config})
are not Boolean but that every row sum is unity.
This result is to be expected
because $\Qmat$ must accommodate interpolation of a constant solution (\eg, $u_{{\rm{g}},i}=1\;\forall i$)
across a nonconforming interface.
\mpar{DLR: discuss the issue of the one common point among 4 elements, of which
one element is refined}

\subsubsection{Global assembly}
\label{sec:global_assembly}

To accommodate Dirichlet boundary conditions \eqref{e:BoCo} into the solution, 
we employ a \emph{masking projection} $\mask$, which is locally diagonal with unit entries everywhere except
corresponding to nodes on Dirichlet boundaries, where there are zero entries.
Any field $\tr{\smash{\vec{\lavec{\phi}}}}\unum=\u\in\uspace_{\vec{b}}$ may be analyzed as
$\u=\u_{\rm{h}}+\u_{\rm{b}}$, where
$\unum_{\rm{h}}\equiv\Jmat\mask\conlab{\Qmat}\unum_{\rm{g}}$ constructs the
projection $\u_{\rm{h}}\equiv\tr{\smash{\vec{\lavec{\phi}}}}\unum_{\rm{h}}\in\uspace_{\vec{0}}$ of $\u$, that is, its
homogeneous part, and $\unum_{\rm{b}}\equiv\unum-\unum_{\rm{h}}$ constructs
$\u_{\rm{b}}\in\uspace_{\vec{0}}$, which vanishes at the interior nodes $\pn[\vec]{\vec{\jmath},k}\not\in\partial\pdomain$ of $\pdomain$.
Inserting this analysis into \eqref{eq:vdiscrete_burgers}
(noting $\v\in\uspace_{\vec{0}}\Rightarrow\vnum=\Jmat\mask\conlab{\Qmat}\vnum_{\rm{g}}$) and repeating the time discretization leading to 
\eq{eq:time_desc_burgers}, we arrive at the following linear equation to solve for $\unum_{\rm{h}}$
at each time step:
\be
\tr{\vnum}\Hmat\unum=\tr{\vnum}\fnum\quad\forall\vnum_{\rm{g}}\Longrightarrow%
\conlab{\tr{\Qmat}}\mask\tr{\Jmat}\Hmat\Jmat\mask\conlab{\Qmat}\unum_{\rm{g}}=%
\conlab{\tr{\Qmat}}\mask\tr{\Jmat}(\fnum-\Hmat\unum_{\rm{b}}),
\label{eq:multi_helm_global}
\ee
where we have ascribed all explicit past-time references from the time-derivative 
expansion and the advection terms in \eq{eq:time_desc_burgers} to $\fnum$.
Equation \eqref{eq:multi_helm_global} is solved by using the PCG algorithm.
While \eqref{eq:multi_helm_global}
shows explicitly that the matrix on the left is symmetric nonnegative-definite\mpar{It is not SPD}, it is not in a form easily solved
on a parallel system.
Left-multiplying \eqref{eq:multi_helm_global} by $\Jmat\mask\conlab{\Qmat}$ and letting
\be
\dss\equiv\Jmat\mask\conlab{\Qmat}\conlab{\tr{\Qmat}}\mask\tr{\Jmat},
\label{eq:dss}
\ee
we get the following local form:
\be
\dss\Hmat\unum_{\rm{h}}=\dss(\fnum-\Hmat\unum_{\rm{b}}).
\label{eq:local_multi_helm}
\ee
The \emph{direct stiffness summation} (DSS) matrix $\dss$
is coded so that the gather and scatter are performed in one operation 
(see \sr{sec:data_comm}), which reduces communication overhead on parallel systems 
\cite{Tufo99}.

In addition to $\mask$ we must introduce the 
inverse \emph{multiplicity matrix} $\mult$ to maintain $\Hone(\pdomain)$ continuity (see \sr{sec:linear_solver}). 
This matrix is diagonal, computed by initializing a collocated vector
$g^\si_{\vec{\jmath},k}=1\;\forall\vec{\jmath},k,\si$, setting child face nodes to $0$,
performing $\lavec{g}\leftarrow\Jmat\conlab{\Qmat}\conlab{\tr{\Qmat}}\tr{\Jmat}\lavec{g}$, then setting
$$
W^{\si,\si'}_{\vec{\jmath},k,\vec{\jmath}\,',k'}=\begin{cases}%
\delta^{\si,\si'}/g_{\vec{\jmath},k}^\si,&
\text{if $\pn[\vec]{\vec{\jmath},k}=\pn[\vec]{\vec{\jmath}\,',k'}$ coincides with a global node,}\\
0,&\text{otherwise.}\end{cases}
$$
For example, corresponding to Figs.\ \ref{fig:conf_config} and \ref{fig:nonconf_config} the diagonals of $\mult$ are
\begin{align*}
&(1,1,\half,1,1,\half,1,1,\half,\half,1,1,\half,1,1,\half,1,1)\\
\text{and}\quad&\textstyle(1,1,1,1,1,1,1,1,1,0,1,1,0,1,1,0,\half,\half,0,\half,\half,0,1,1,0, 1,1),
\end{align*}
respectively.

\subsection{Linear solver}
\label{sec:linear_solver}

The modifications to the PCG algorithm required to solve
\eq{eq:local_multi_helm} in the nonconforming case stem
from the requirement that the iteration residuals $\rnum$
and the search directions $\wnum$ correspond to functions
$\vec{r}\equiv\tr{\smash{\vec{\lavec{\phi}}}}\rnum$ and
$\vec{w}\equiv\tr{\smash{\vec{\lavec{\phi}}}}\wnum$ belonging to
$\Hone(\pdomain)^d$.
The CG algorithm searches the global \dofs\ space for the solution
to the linear equation.
So that we may continue to use the local matrix forms, however, we must
also mask off all Dirichlet nodes (if any exist), which are
not solved for.
The assembly matrix masks off these nodes in such a way that the
new search direction $\vec{w}\in\Hone(\pdomain)^d$.
Additionally, in all cases in the CG iteration where a quantity
$\vec{g}$ must remain in $\Hone(\pdomain)^d$, we explicitly ``smooth''
it by $\lavec{g}\leftarrow\smooth\lavec{g}$, using an $\Hone$ smoothing matrix
$\smooth\equiv\Jmat\mask\conlab{\Qmat}\conlab{\tr{\Qmat}}\mult$.
The result of the operation is a quantity that is interpolated
properly to the child edges and that is expressed precisely
once (unlike in the DSS case) at multiple local nodes that
represent the same spatial location.
A similar smoothing matrix
$\smooth_{\rm{f}}\equiv\Jmat\conlab{\Qmat}\conlab{\tr{\Qmat}}\mult$,
where the Dirichlet mask is not inserted, is used to smooth
the final inhomogeneous solution.
Note that it is critical that the inhomogeneous boundary term
$\u_{\rm{b}}\in\Hone(\pdomain)^d$ in \eq{eq:local_multi_helm};
thus, the smoothing matrix $\smooth$ is applied to $\unum_{\rm{b}}$
before the Helmholtz operator is.
However, the nonsmoothed boundary term must be added after the
convergence loop in order to complete the solution.
Note also that the final smoothing operation follows the addition
of the boundary condition and therefore \emph{cannot} be masked;
hence the distinction of the final $\smooth_{\rm{f}}$ matrix.

\subsubsection{Modified preconditioned conjugate-gradient algorithm}
\label{sec:pcgalgorithm}

\begin{figure}
\renewcommand{\arraystretch}{.7}
\begin{tabular}{ll}
$\unum_{\rm{h}}=\boldsymbol{0}$&
{\tt\cmt\ initialize homogeneous term}\\
$\rnum = \dss(\fnum - \Hmat\smooth\unum_{\rm{b}})$&
{\tt\cmt\ initialize residual}\\
$\wnum=\boldsymbol{0}$&
{\tt\cmt\ initialize search vector}\\
$\rho_1=1$&
{\tt\cmt\ initialize parameter}\\
\multicolumn{2}{l}{\tt Loop until convergence:}\\
$\qquad\enum = \smooth\Pmat^{-1}\rnum$&
{\tt\cmt\ error estimate}\\
$\qquad\rho_0 = \rho_1$, $\rho_1 = \tr{\rnum}\mult\enum$&
{\tt\cmt\ update parameters}\\
$\qquad\wnum\leftarrow\enum + \wnum\rho_1/\rho_0$&
{\tt\cmt\ increment search vector}\\
$\qquad\rnum' = \dss\Hmat\wnum$&
{\tt\cmt\ image of $\wnum$}\\
$\qquad\alpha = \rho_1/\tr{\wnum}\mult\rnum'$&
{\tt\cmt\ component of $\unum_{\rm{h}}$ increment}\\
$\qquad\unum_{\rm{h}}\leftarrow\unum_{\rm{h}} + \alpha \wnum$&
{\tt\cmt\ increment $\unum_{\rm{h}}$ along $\wnum$}\\
$\qquad\rnum\leftarrow\rnum - \alpha \rnum'$&
{\tt\cmt\ increment residual}\\
\multicolumn{2}{l}{\tt End Loop}\\
\multicolumn{2}{l}{$\unum = \smooth_{\rm{f}}(\unum_{\rm{h}} + \unum_{\rm{b}})$.}
\end{tabular}
\caption{PCG algorithm}
\label{t:pcgalgorithm}
\end{figure}
%
With these considerations we present in Table
\ref{t:pcgalgorithm} the PCG algorithm for
the assembled local problem \eq{eq:local_multi_helm} for conforming elements
\cite[see][Sec. 4.5.4]{DFM2002} modified for nonconforming
elements.
Preconditioning is handled by the matrix $\Pmat^{-1}$.
In general, the preconditioned quantity must be smoothed.

\subsubsection{Preconditioners}
\label{sec:preconditioners}
Much investigation remains to determine the optimal preconditioners
$\Pmat^{-1}$ for the equation solvers supported in \gaspar.
However, the preconditioning strategy is not the present work's
focus.
Nevertheless, for the purpose of turbulence study we, make some
general comments to guide the development of optimal preconditioners.
If we consider the expression \eq{eq:helmholtzop} for the
spectral-element Helmholtz operator, we note that $\nu\propto\Rn^{-1}$
is very small.
Also, because we always use relatively high spatial degree, the
timestep $\varDelta{t}\propto1/\beta_{\rm bdf}^{n,n}$ is very small.
Hence, we see that $\Hmat_k^n$ is strongly diagonally dominant.
Thus, approximating $\Hmat_k^n$ with its diagonal entries alone is a
reasonably good choice for a preconditioner, and this is used in
the present work.

\subsection{Adaptive mesh formulation}
\label{sec:adaptivity}
As mentioned in \sr{sec:variational_form}, the global domain
$\pdomain$ is initially covered \eqref{e:domdecom} by a set of disjoint
(nonoverlapping) elements $\elbox{k}$.
Each of these initial elements becomes a tree root element, which
is identified by a unique root \emph{key} $k_{\rm{r}}\isin{1,2,3}{}$
for that tree.
At each level $\ell\isin{\ell_{\rm{min}}}{\ell_{\rm{max}}}$, an element data
structure provides both its own key $k$ and its root key $k_{\rm{r}}$.
For any level $\ell$, the range of $2^{d\ell}$ valid element keys will
be $k\in[2^{d\ell}k_{\rm{r}},2^{d\ell}(k_{\rm{r}}+1)-1]$ because the
refinement is \emph{isotropic} (that is, it splits an element at the
midpoints of all its edges to produce its $2^d$ child elements).
Conversely, we obtain the level index from the element key using
\begin{equation}
\ell=\lfloor\log_{2^d}(k/k_{\rm{r}})\rfloor.
\label{e:keytolev}
\end{equation}
In order to ensure all keys are unique, given a root $k_{\rm{r}}$ the next root is
$k_{\rm{r}}'=2^{d\ell_{\rm{max}}}(k_{\rm{r}}+1)$, and so on.

After elements $\elbox{k}$ are identified (``tagged'') 
for refinement or coarsening at level $\ell$, three
steps are involved in performing DARe:
(1) performing refinement by
adding a new level of $2^d$ child elements $\elbox{2^dk},\cdots\elbox{2^d(k+1)-1}$ at level $\ell+1$ to replace
each $\elbox{k}$, or else coarsening $2^d$
existing children $\elbox{k},\cdots\elbox{k+2^d-1}$ into a new parent $\elbox{\lfloor{k}/2^d\rfloor}$; (2) building data structures for
all element boundaries, which hold data representing global \dofs\ and accept
gathers ($\tr{\Qmat}\unum$ values) or perform scatters
($\Qmat\unum_{\rm{g}}$ values); and (3) determining neighbor lists for data exchange.
Neighbor lists consist of records (structures) that each contain
the computer processor id, element key $k$, root key  $k_{\rm{r}}$ and boundary id $s\isin{0}{2^d-1}$ of each neighbor element that
adjoins every interface. 
In refining or coarsening, the field
values for each child (parent) elements are interpolated from the parent (child) fields.
For simplicity, the interior of each element boundary is restricted to
an interface between one coarse and at most $2^{d-1}$ refined neighbors. Thus,
at most one refinement-level difference will exist across the interior of
an interface between neighboring elements.

In \gaspar, the data structures that represent global \dofs\ at the inter-element
interfaces are referred to as ``mortars.'' These structures are not to
be confused with the mortars used in MEM; however, they
serve as templates for that more general method.
Recalling \fig{fig:nonconf_config} as a paradigm, in general the
mortars contain node locations and the basis set of
the parent edge (in 2D, or face in 3D). These structures represent the same
field information
for the parent and child edges; their nodes coincide with the parent edge's nodes, and they
interpolate global \dof\ data to the child edges,
as described above. The mortar data structures are determined by communicating
with all neighbors to determine which interfaces are nonconforming. This
communication uses a \emph{voxel database} (VDB)
\cite{henderson95}. A voxel database consists of records containing
geometric point locations, a component id that tells what part of 
the element $\elbox{k}$ (edge $\partial\elbox{k,s}$, vertex $\in\partial^2\elbox{k,s}$, etc.)
the point represents, an id of the element that contains the point,
that element's root id, and
some auxiliary data. Two VDBs are constructed: one consists of all element
vertices, and one consists of all element edge midpoints. With these
two VDBs, we are able to determine whether a relationship between
neighbor edges is conforming and also determine the mortar's geometrical
extent. The VDB approach can also be used for general deformed geometries
in two and three dimensions, as long as adjacent elements share well-defined
common node points. 

The algorithm classes that carry out refinement operate only on the element
and field lists. The SEM solvers adjust themselves automatically to accommodate
the lists that are modified as a result of DARe.

\subsubsection{Refinement and coarsening rules}
\label{sec:refinerules}
The method that actually carries out the
refinement and coarsening takes as arguments only two buffers: one containing the local
indexes of the elements to be refined, and one with indexes of elements to be coarsened.
While the refinement criteria that identify elements for
refinement or coarsening are described in \sr{sec:error},
we point out here that before mesh refinement or coarsening is undertaken, 
the tagged elements are checked for compliance with several rules. For refinement,
we have the following rules:
\begin{enumerate} 
\item[r1.] The refinement level must not exceed a specified value.
\item[r2.] No more than one refinement level may separate neighbor elements.
\end{enumerate} 
These rules must be adhered to even in the case of interfaces that lie on periodic boundaries.
The refinement list is modified to remove any elements that would violate 
rule r1, and rule r2 is enforced by tagging a coarse
element for refinement if it has a refined neighbor that is in the current
refinement list. This is most easily effected by building a global
refinement list, consisting of the element keys of all elements tagged for
refinement.
Each local element's neighbors are then checked; if a neighbor
is in the global refinement list and it is a nonconforming neighbor,
then the local element must also be refined.

We may not coarsen an element if under any of the following conditions:
\begin{enumerate} 
\item[c1.] The element is the tree root.
\item[c2.] Any of its $2^d-1$ siblings are not tagged for coarsening.
\item[c3.] The element appears in a refinement list.
\item[c4.] Refinement rule r2 would be violated.
\end{enumerate} 
To enforce rule c4, we introduce the 
notion of a \emph{query-list}, in other words, a global list of records
of each element key $k$, its parent key $\lfloor{k}/2^d\rfloor$,
and its refinement level $\ell$ \eqref{e:keytolev}.
The global list is built from local buffers of element indexes that are 
then gathered from among all processors. 
The following procedure
is then used: 
\begin{enumerate}
\item[1.]
Build a query-list (RQL) from the element indexes in the refinement list.
\item[2.]
 Find the maximum and minimum refinement levels $\ell_{\rm{max}}$ and $\ell_{\rm{min}}$ represented among
all keys tagged for coarsening. 
\item[3.]
Reorder the current local coarsen list from $\ell_{\rm{max}}$ down to $\ell_{\rm{min}}$.
\item[4.]
Working from $\ell=\ell_{\rm{max}}$ down to $\ell_{\rm{min}}$, 
\renewcommand{\labelenumii}{\roman{enumii}.}
\begin{enumerate}
\item Build a query-list (CQL) from the element indexes in the current coarsen list.
\item For all elements in the local coarsen list at the current $\ell$, check that (a) any 
refined neighbor is in the CQL and (b) no refined neighbors are in the the RQL.
If both conditions are met at this $\ell$, then the local element
index is retained in the current coarsen list. Otherwise it is deleted.
\end{enumerate}
\item[5.]
Do a final check that all elements in the local coarsen list have all their siblings also tagged for
coarsening. This is done by building a query-list of all current coarsen lists and verifying 
that a local element's siblings are in the global list. Sibling elements are identified by
having the same parent key as the local element.
\end{enumerate}

We note that in order to make the algorithm consistent, the local refinement lists are
checked, and possibly modified \emph{before} checking and modifying the coarsen lists.

\subsubsection{Communicating boundary data}
\label{sec:data_comm}

The mortar data structures contain all the data to be communicated
between elements during each application of the DSS or smoothing operation \eqref{eq:dss}. 
\gaspar\ communicates element-boundary data by exchanging data between adjacent elements, 
which necessitates network communication on parallel computers. This data exchange
involves an {\it initialization} step and and {\it operation} step. The
initialization step establishes the required element/processor connectivity 
by performing a bin-sort of global node indexes and having each processor process
the nodes from a given bin to determine neighbor lists. This method was 
suggested in \cite[ch.\ 8]{DFM2002} but to our knowledge has never 
been implemented. The procedure
is to label the mortar-structure nodes with unique
indexes, generated by computing the Morton-ordered index for each geometric node point.
This index is computed by integralizing the physical-space coordinates and then interleaving
the bits of each integer component to create a unique integer.  Thus, all nodes representing the same geometric position will have the same index label.
For $P$ processors, a collection of bins $\set{B}_l$, $l\isin{0}{P-1}$, is generated that
partitions the dynamic range (global maximum)\mpar{What does ``global maximum'' mean here?} of the node indexes.
Each processor $l$ partitions its list of node indexes into the bins, sending
the contents of bin $\set{B}_{l'}$ to processor $l'$, where the data are combined with those
from other processors and then sent back to the originating processor. After this
step, a given processor is informed of which other processors share which nodes among all
the mortar nodes. 

With the information gleaned from this initialization step,
the operation step involves the communication
of the data at any node point with all other processors that share
that node. This data is extracted from the element in question by using the
pointer indirection provided by the local-to-global map represented by the unique
node index. The data at common vertices is summed for the DSS or smoothing
operations and returned to the local node also by indirection. The algorithm
provides that common vertices residing on the same processor are summed
before being transmitted to the other processors that share the node, in order 
to reduce the amount of data being communicated. At the end of
the operation step, the data at multiply-represented global nodes are identical.
This gather-scatter procedure
ensures that the DSS'ed data are available for local computation immediately after
the final communication of data.

One benefit of this method for performing the gather-scatter operations is that it allows
the communication to separated from the geometry because the global unique node ids
are essentially unstructured lists of local data locations. However, the method is
somewhat inefficient in the current formulation of DARe in \gaspar.  
Neighbor lists for transmitting data may also be constructed 
by using the VDB. Since the VDBs are synchronized and thus represent, in a sense, global
data, this same VDB can also be used to determine a given
element edge's neighbor lists. After the mortars are constructed, each element
edge is associated with its neighbors' local ids and processor ids. Hence, the 
neighbor lists for handling element-boundary data exchange are determined easily 
from existing global data; there is no need to perform the initialization step
in the gather-scatter method described above because the information is readily
available from the VDB synchronization.

\subsubsection{Error estimators}
\label{sec:error}
Elements are tagged for refinement or coarsening by using an a posteriori refinement
criterion. One criterion, adapted from \cite{mavriplis93}, 
uses the Legendre spectral information contained within the spectral-element representation to
estimate the quadrature and truncation errors in the solution and to compute the rate at which
the solution converges spectrally in each element $\elbox[\bar]{k}$. We refer to this as the 
{\it spectral estimator}. The discretized solution along 1D lines in coordinate
direction $\si'$ (fixing the other coordinates) is represented as a Legendre expansion
\be
u^\si(\comap[\vec]{k}(\cdots,\xi^{\si'}\cdots))=\sum_{j=0}^{\pode}\check{u}_j^\si\Leg{j}(\xi^{\si'}),
\label{eq:legendre_exp}
\ee 
where the $\check{u}_j^\si$ can be computed easily by using the orthogonality of $\Leg{j}$ 
w.r.t.\ \eqref{e:innpro}.
Then an approximation to the solution error along the line can be expressed as
\be
\varepsilon_{\rm est}^\si=\sqrt{\sum_{j=\pode}^\infty\left(j+\half\right)^{-1}(\check{u}_j^\si)^2},
\label{eq:spect_error}
\ee
where the $j=\pode$ term is the (over)estimate of the quadrature error and
the $j>\pode$ terms sum to 
the truncation error.  Since we do not have $\check{u}_{j>\pode}^\si$, we must extrapolate
using the coefficients we do have. 
We assume $\check{u}_{j>\pode}^\si$ can be approximated by
$\ln|\check{u}_j^\si|\approx\ln C^\si-\lambda^{\si}j$, and we determine $C^\si$ and $\lambda^{\si}$ from a
linear fit using the $M$ final
coefficients $\check{u}_{\pode-M<j\leq\pode}^\si$ in \eqref{eq:legendre_exp}. With this continuous approximation 
for the extrapolated coefficients, we approximate the truncation error term in
\eq{eq:spect_error} as an integral over $j$ and integrate analytically to compute the total error $\varepsilon_{\rm est}^\si$.
In this work, the maximum
error over all $\nopo$ lines for each $\si'$ within $\elbox[\bar]{k}$ is 
computed and taken as the $\elbox[\bar]{k}$-local error estimate. The local convergence rate is
simply the minimum value of $|\lambda^{\si}|$ over all lines and all $\si'$. For all our test cases, $M=5$.

Thus, $\elbox[\bar]{k}$ is refined, if for some $\si$, $\varepsilon_{\rm est}^\si$ is above a threshold value or if $|\lambda^{\si}|$ is below
another threshold. For coarsening, for all $\si$, all $2^d$ sibling elements must have their $\varepsilon_{\rm est}^\si$s below some value
proportional to the refinement threshold. This requirement prevents ``blinking,''
where refined elements are immediately coarsened because the error tolerances
are met on the refined elements.

In conjunction with this spectral criterion, 
we can often obtain better overall accuracy-convergence results by
checking whether the $\elbox[\bar]{k}$-maximum second derivative in any coordinate is above a certain threshold
and by performing a logical {\rm OR} of that condition with the spectral 
error threshold for refinement tagging. 

While the high degrees used in the expansions
will help the spectral error estimator,
other refinement criteria may be more effective, given the variety of solution structures 
arising in our applications.
The investigation of refinement criteria appropriate for such intermittent features 
is a major outstanding problem in numerical solution of PDEs, and one that we consider
in a companion paper \cite{FRP2005}.

\section{Test problems and results}
\label{sec:test_problems}

We have chosen a set of test problems that examine various aspects of \eq{eq:burgers}.
The primary objective is to investigate the temporal and spatial convergence
properties of the solutions when adaptivity is used. For this purpose, we have
selected test problems that have analytic solutions, so that the errors may
be determined exactly, instead of only by comparison, for example, to a highly refined
control solution. The test problems begin with the simplest aspect of
\eq{eq:burgers} and continue to progressively more difficult problems
until the behavior of the full 2D nonlinear version of \eq{eq:burgers} is considered.

For each test problem, a BDF3/Ext3 scheme is used for the time derivative and the 
advection term, respectively (\sr{sec:time_discrete}).  This requires that, at the initial time
$t^0$, all the 
required time levels $t^m$ be initialized, $m\isin{1}{\max(M_{\rm{bdf}},M_{\rm{ext}})-1}$.
Both the spectral and second-derivative 
error estimators are used for adaption criteria. The spectral error
is normalized by the norm of the solution, $\dlebnorm{\unum^0}\equiv\sqrt{\tr{{\unum^0}}\unum^0}$,  at the 
start of the run, and the second $x^\si$-derivative is normalized by 
$\dlebnorm{\unum^0}/L^2$, where $L$ is the longest global domain length. Elements
are tagged for refinement based on a logical OR of the two criteria.
The $|\lambda^{\si}|$ threshold in all cases is $\ln10$.

\subsection{Heat equation}
\label{sec:heateqtest}

For the linear case $\ua=\vec{c}(t)$ the analytic fundamental solution of \eq{eq:burgers} is
a $d$-periodized Gaussian in $\pdomain=[0,1]^d$:
\be
u^\si_{\rm{a}}(\x,t)\equiv\frac{\sigma(0)^d}{\sigma(t)^d}\sum_{\imath^1,\cdots\imath^d=-\infty}^{\infty}%
\exp-\left(\frac{\x-\x^{\,0}+\vec{\imath}-\int_0^t\vec{c}(t')\d{t'}}{\sigma(t)}\right)^2
\label{eq:sph_gaussian_advect}
\ee
for $t>-\sigma(0)^2/4\nu$ ($u^\si_{\rm{a}}(\x,t)\equiv0$ otherwise),
where $\sigma(t)\equiv\sqrt{\sigma(0)^2+4\nu t}$,
$\sigma(0)=\sqrt{2}/20$ is the initial distribution e-folding width
and $\x^{\,0}=\sum_{\si=1}^d\uv{\si}/2$ is the initial peak location.
To compute \eqref{eq:sph_gaussian_advect}, we truncate
summands of value be less than $10^{-18}$ of the partial sum.
The simplest version of \eq{eq:burgers} is the heat equation, where
$\ua=\vec{0}$. The purpose here is to investigate temporal and spatial convergence
of the adaptive solutions without advection.
The initial condition at $t=0$ from \eqref{eq:sph_gaussian_advect} with $d=2$
is computed on $K=4\times4$ elements, and the mesh is refined until the
maximum allowed number of refinement levels may be reached. Both a spectral estimator
and derivative threshold were used. The threshold and coarsening 
factor for each were set to be $10^{-4}$ ($10^{-2}$) and $1$ ($0.5$), respectively.

\subsubsection{Temporal convergence}
\label{sec:heateqtest:tconv}
\begin{figure}
\begin{center}
\includegraphics[scale=.60]{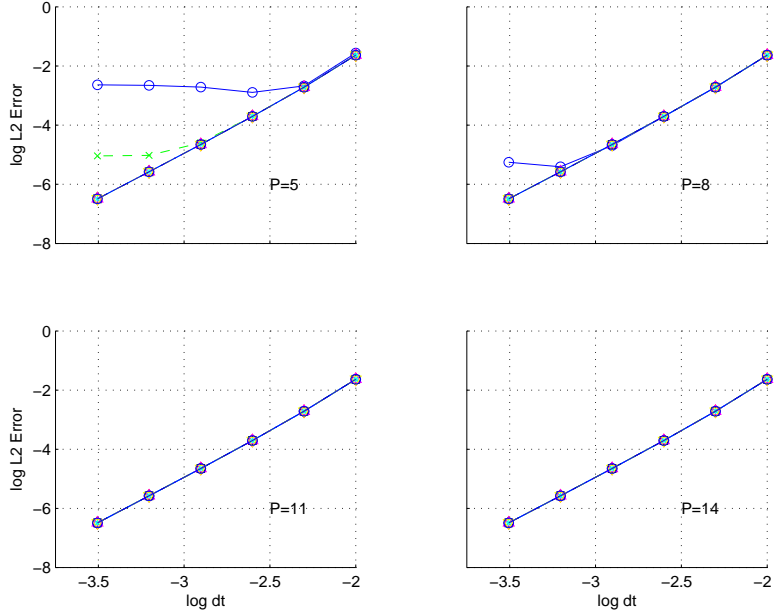}
\end{center}
\caption{Plots of $\log_{10}(\dlebnorm{\unum-\unum_{\rm{a}}}/\dlebnorm{\unum_{\rm{a}}^0})$ for the heat equation
vs $\log_{10}\varDelta t$, for different polynomial degrees $\pode$.
Within each plot are four curves for the maximum refinement levels $\ell_{\rm{max}}=0$ (solid curve with 
circle markers), $\ell_{\rm{max}}=1$ (dashed curve with cross markers), $\ell_{\rm{max}}=2$ (dotted curve with 
diamond markers), and $\ell_{\rm{max}}=3$ (dash-dotted curve with square markers). 
The control solutions are indicated with dashed curves and follow closely the
$\ell_{\rm{max}}=3$ curves. The axes in each plot have the same limits. Note that as $\pode$ increases, 
the curves converge, for this range of $\varDelta t$.
\label{fig:heat_dtconv}
\bigskip
}
\end{figure}

First, we consider time convergence of the adaptive solutions by integrating to 
a fixed time $t=0.05$ for various fixed timesteps $\varDelta t$. From \eq{eq:sph_gaussian_advect} a curve of relative $\Ltwo$ error
$\varepsilon=\dlebnorm{\unum- \unum_{\rm{a}}}/\dlebnorm{\unum_{\rm{a}}^0}$
vs $\varDelta t$ is generated for each of several maximum-refinement levels $\ell_{\rm{max}}$ and 
for four degrees $\pode$. We present the results in \fig{fig:heat_dtconv}, 
adopting the convention $\ell$-\emph{control} for the grid that uniformly covers the domain with 
elements at the finest resolution in the $\ell_{\rm{max}}=\ell$ case.
The BDF3/Ext3 is a globally third-order scheme, so we expect that if the solution is well resolved
spatially, we should see a slope of $\approx 3$ in a log-log plot of error vs $\varDelta t$.
This is precisely what is seen in the figure; each plot consists of a sequence of 
four curves for the refinement levels $\ell_{\rm{max}}\isin{0}{3}$, where $\ell_{\rm{max}}=0$ implies that
no refinement is done.
For the spatially resolved curves in each plot, the error is linear with slope 3.14. 

We see in \fig{fig:heat_dtconv} that even at low $\pode$, the solution can be well resolved as 
long as refinement is used. We also see that as $\pode$ increases, there is less need for
refinement because the unrefined mesh is able to resolve the solution adequately, at least 
for the period of integration. 

To compare the various refinement levels, we have made the same runs using a
3-control grid.
The $\ell_{\rm{max}}$-control solutions are defined to be those generated on a nonadaptive
grid comprising elements at the finest scale in the $\ell_{\rm{max}}$-adaptive
case. Thus, in this problem, since we initialized with a $K=4\times4$ grid, the
$3$-control grid consists of $2^3 4 \times 2^3 4$ or $K=32\times32$ elements. 
These plots are indicated
by the dashed lines in \fig{fig:heat_dtconv},
which all follow the $\ell_{\rm{max}}=3$ curves. As the polynomial degree increases, the 
solution requires fewer levels of refinement to achieve the same accuracy as
with the 3-control grid, at a considerable computational savings. 

\subsubsection{Spatial convergence}
\label{sec:heateqtest:spconv}
\begin{figure}
\begin{center}
\includegraphics[scale=.60]{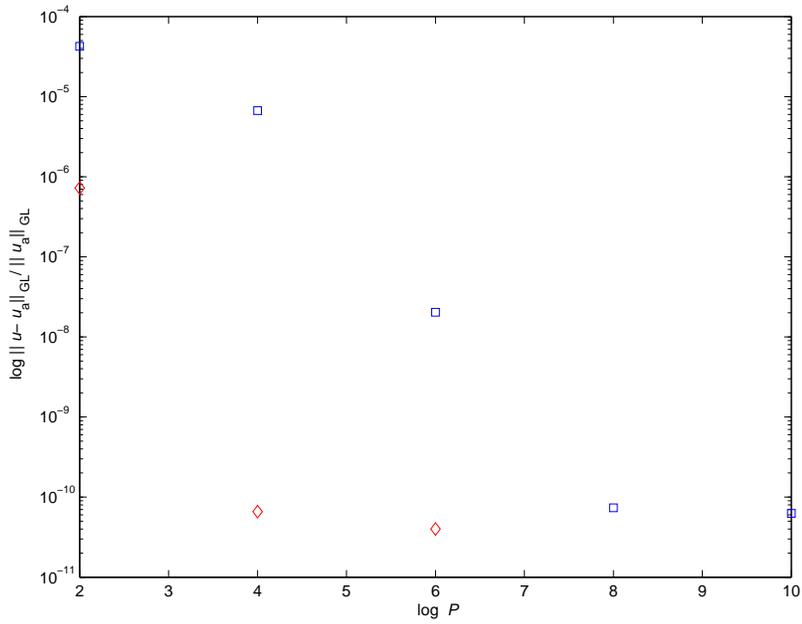}
\end{center}
\caption{Plots of $\log_{10}(\dlebnorm{\unum-\unum_{\rm{a}}}/\dlebnorm{\unum_{\rm{a}}^0})$, as a function
of $\pode$ for the heat equation test. The square markers indicate the adaptive runs,
 while the diamonds represent the 
3-control grid runs at the same polynomial degree. 
\label{fig:heat_spconv}
}
\end{figure}

In this portion of the heat-equation test, we consider the effects of 
polynomial degree $\pode$ on the solution. The maximum number of refinement levels is fixed
to $\ell_{\rm{max}}=3$. Here, a variable Courant-limited timestep
$$
\varDelta{t}\leq\kappa\left/\max_{j\isin{1}{\pode},k\isin{1}{K},\si\isin{1}{d}}\left(%
\frac{4\nu}{(\elsiz[\si]{j,k})^2}+%
\frac{|u_{j-1,k}^\si|+|u_{j,k}^\si|}{2\elsiz[\si]{j,k}}\right)\right.
$$ is used with a
fixed Courant number of $\kappa=0.2$,
where $\elsiz[\si]{j,k}\equiv|\pn{j,k}^\si-\pn{j-1,k}^\si|$ uses \eqref{e:npmapped}. The solution is integrated to a time $t=0.5$ chosen so that 
we observe the solution coarsening as it decays.
The initial mesh is the same as in the time convergence
test. The spatial convergence result is presented in \fig{fig:heat_spconv}.

\begin{figure}
\begin{center}
\includegraphics[scale=.60]{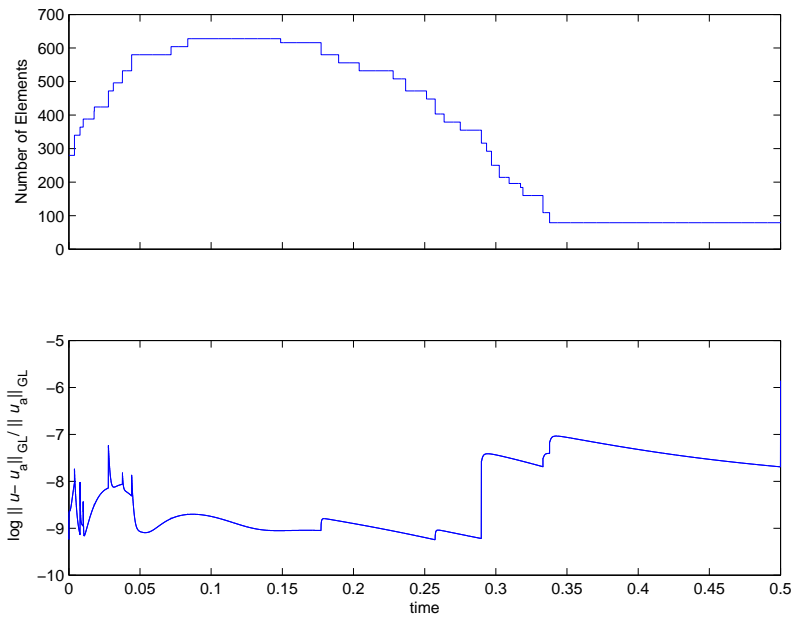}
\end{center}
\caption{Plots of diagnostic quantities vs time in the $\pode=6$ heat equation test.{\it Top}: $Number\; of\; Elements$ vs $time$. {\it Bottom}: $\log_{10}(\dlebnorm{\unum-\unum_{\rm{a}}}/\dlebnorm{\unum_{\rm{a}}^0})$ vs $time$. 
\label{fig:heat_tplots}
}
\end{figure}

We see in the figure the linear behavior characteristic of the spatial
convergence in all our test problems.
We expect that a solution converging spectrally as in \eqref{e:errorbound} will exhibit a straight line in a 
plot of $\log_{10}(\dlebnorm{\unum-\unum_{\rm{a}}}/\dlebnorm{\unum_{\rm{a}}^0})$ vs $\pode$,
indicating that the error decays exponentially.
Also plotted in this figure are the 3-control
runs, which saturate very quickly in this case but still provide uniformly
better errors than in the adaptive case except at the highest degrees. The cause is
likely the elliptic nature of the problem in which the coarsening 
elements transmit their error throughout the grid. In \fig{fig:heat_tplots} we
can see that 
the error over time does not decay monotonically except during periods where
the grid is quasi-static. 
Keeping in mind that the
decaying Gaussian contains all polynomial orders, one concludes that the solution 
error is globally influenced by the error from 
interpolations and from the inability of the truncated polynomial 
expansion at this degree to accurately model the
decaying Gaussian.

\subsection{Linear advection equation}
\label{sec:linadvecttest}
In our next test we consider the linear advection equation \eq{eq:burgers} with $d=2$ and constant
$\ua=\uv{1}$. This problem considers the ability of the code to simulate and follow
a reasonably sharply localized translating distribution. 
The initial distribution is given by \eq{eq:sph_gaussian_advect} at $t=0$.
For this test problem, we set $\nu = 10^{-4}$.
The spectral error tolerance in this problem is turned off.
The derivative
criterion is set to $1.0$ with a coarsening factor of $0.5$.
The $|\lambda^{\si}|$ threshold is $\ln10$.

\subsubsection{Temporal convergence}
\begin{figure}
\begin{center}
\includegraphics[scale=.60]{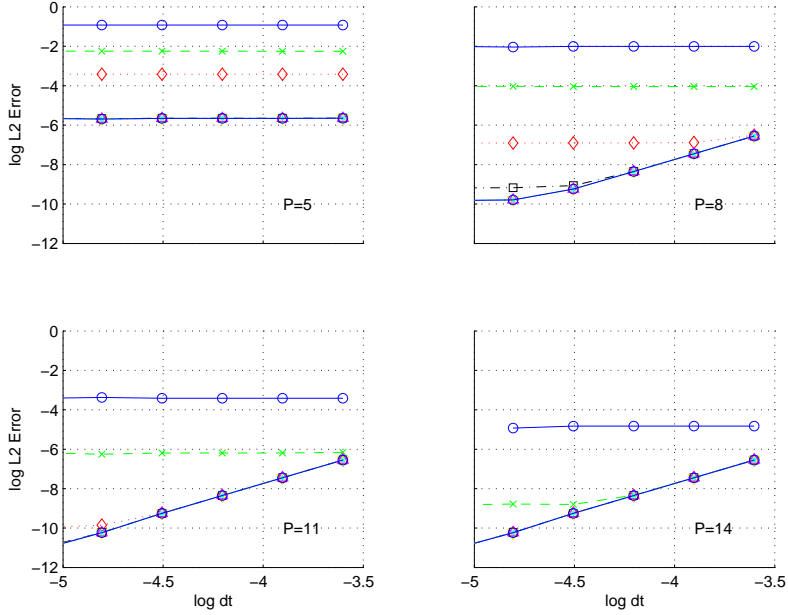}
\end{center}
\caption{Plots of $\log_{10}(\dlebnorm{\unum-\unum_{\rm{a}}}/\dlebnorm{\unum_{\rm{a}}^0})$ for the advection-dominated problem vs $\log_{10}\varDelta t$, for different $\pode$.
Within each plot are subplots for each of four refinement levels, $\ell_{\rm{max}}=0$ (solid curve with circle markers),
$\ell_{\rm{max}}=1$ (dashed curve with $\times$ markers), $\ell_{\rm{max}}=2$ (dotted curve with diamond markers), and $\ell_{\rm{max}}=3$ (dash-dotted curve with square markers). The control solutions are indicated with dashed curves and generally follow the
$\ell_{\rm{max}}=3$ curves. The axes in each plot have the same limits.
\label{fig:galump_dtconv}
}
\end{figure}

The temporal convergence test is done in the same way as for the heat equation
(\sr{sec:heateqtest:tconv}). The total integration time is $t=0.06$.
We begin with a $K=4\times4$ element mesh, each element of 2D degree $\pode\times\pode$. In \fig{fig:galump_dtconv}
we present our results.

For the spatially resolved curves in each plot, the slope of the curve is 2.95. 
Even at high degree $\pode$, the error is $\varDelta t$-independent for the unrefined mesh. 
For lower degrees, to obtain a case where the solution error decays at the order of the 
time-stepping method indeed requires a larger number of refinement levels, indicating that the
solution is well resolved spatially only at these higher levels. Thus, in order to resolve this
distribution properly so that the temporal error is $\oforder{\varDelta t^3}$, refinement is
necessary. 

Also plotted \fig{fig:galump_dtconv} are the 3-control runs corresponding 
to the adaptive solutions.  These runs are indicated
by the dashed lines (just visible in the top-right plot), 
which all follow the $\ell_{\rm{max}}=3$ curves. As the polynomial degree increases, the 
solution requires fewer levels of refinement to achieve the same accuracy as
with the 3-control grid, again at a considerable computational savings. 

\subsubsection{Spatial convergence}
\begin{figure}
\begin{center}
\includegraphics[scale=.60]{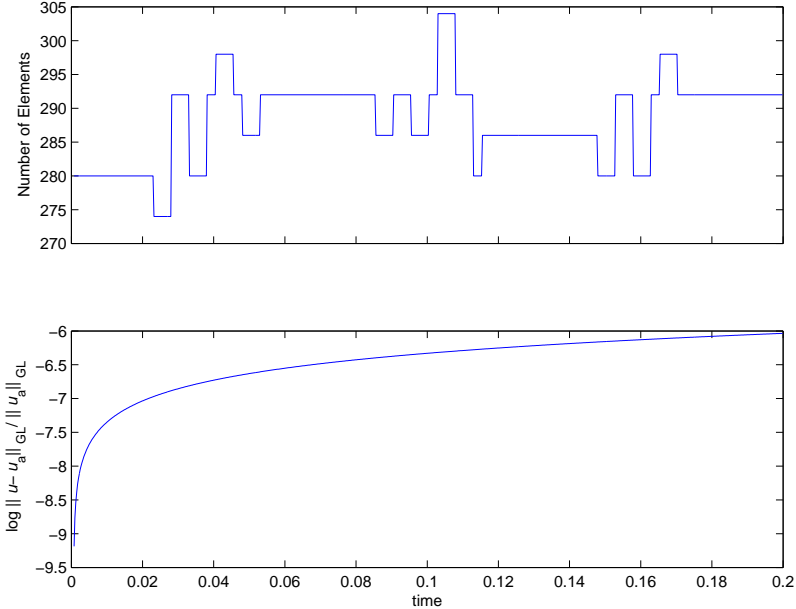}
\end{center}
\caption{Plots of diagnostic quantities vs time in the $\pode=8$ linear advection test.{\it Top}: {\rm number of elements vs time}. {\it Bottom}: $\log_{10}(\dlebnorm{\unum-\unum_{\rm{a}}}/\dlebnorm{\unum_{\rm{a}}^0})$ vs $time$. The errors for the adaptive and control meshes lie on top of one another. 
\label{fig:galump_tplots}
}
\end{figure}

We next consider the effects of expansion degree $\pode$ on the solution error.
The maximum number of refinement levels is fixed
to $\ell_{\rm{max}}=3$. Here, a Courant-limited timestep is again used with a
Courant number of $0.2$. The solution is integrated to a time ($t=0.2$) chosen so that several
cycles of coarsening and refinement occur (see top of \fig{fig:galump_tplots}). The initial mesh is the same as in the time convergence
test. The spatial convergence result is presented in \fig{fig:galump_dtconv}.
\begin{figure}
\begin{center}
\includegraphics[scale=.60]{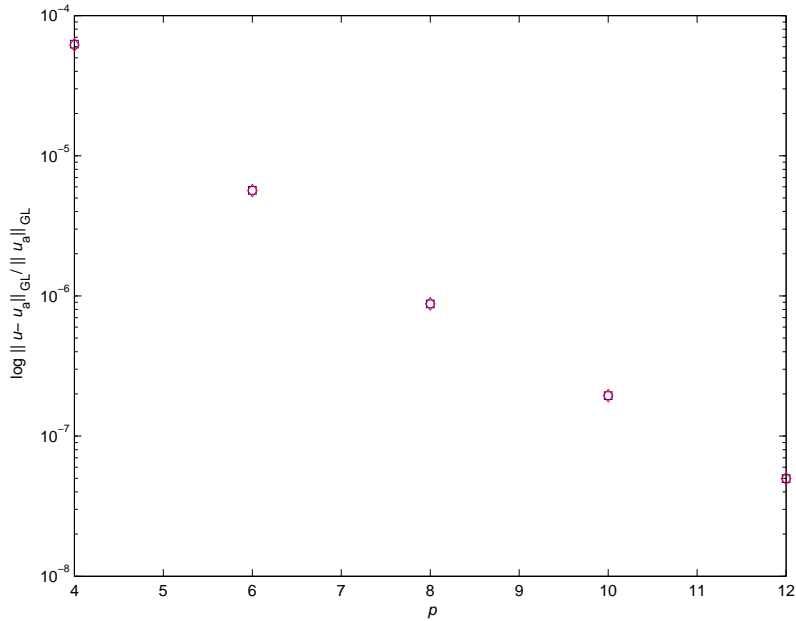}
\end{center}
\caption{Plots of $\log_{10}(\dlebnorm{\unum-\unum_{\rm{a}}}/\dlebnorm{\unum_{\rm{a}}^0})$ as a function
of $\pode$ for the linear advection test. The square markers indicate the adaptive runs, 
while the diamonds represent the 
3-control grid runs at the same polynomial degree. These two curves are nearly identical. 
\label{fig:galump_spconv}
}
\end{figure}

The anticipated spectral decay of error can be seen in \fig{fig:galump_spconv}, which also 
includes the uniform control solutions. 
We see again that the adaptive solution error decays nearly identically 
to the control solution, suggesting again that interpolations and 
polynomial truncation error introduce
no deleterious effects in the linear advection case.

In \fig{fig:galump_tplots}, typical plots of two diagnostics are provided for the $\pode=8$ run in 
\fig{fig:galump_spconv}. In the top of the figure is shown the number of elements as a function 
of time, and in the bottom, the error. Clearly, the adaption does not alter the monotonic
behavior of the error. The error for the 3-control grid in the $\pode=8$ case is
also plotted in \fig{fig:galump_tplots}; the behavior of this error as a function
of time is nearly identical to that for the adaptive mesh. Since the problem was initialized 
with a $K=4\times4$ grid at $\ell_{\rm{min}}=0$, at $\ell_{\rm{max}} = 3$, the 3-control problem has $K=32\times32$ elements. The 
adaptive case clearly provides for a significant savings in the number of
\dofs\ required to compute the same solution. 

\subsection{2D Burgers equation}
\label{sec:burgerstest}

In this section we examine the full nonlinear $\ua=\u$ version of \eq{eq:burgers}.  
The objective is to investigate the 
solution errors as the mesh resolves or tracks stationary and propagating fronts
that are generated and sustained by the nonlinearity of the equation.
The first problem is the familiar Burgers stationary front, and the second,
a radial N-wave. With the former, we investigate the effects of 2D
adaptivity on a rotated 1D nonlinear problem. In the N-wave
problem we consider convergence properties of a fully 2D nonlinear case.

\subsubsection{Stationary Burgers front problem}

The stationary Burgers front problem is the classical solution to the nonlinear 
advection-diffusion equation
\eq{eq:burgers} with a planar front developing in the $x$ direction. We compare the 
maximum derivative of the field, $|\PD{u}{x}|_{\rm max}$, and the time at which the maximum
occurs with the analytic solution.
The classical 1D Burgers front problem for $q(y,t)$ is
cast into a 2D framework by the substitution
$\u(\x,t)=\vec{k}q(\DP{\vec{k}}{\x},k^2t)$ in \eq{eq:burgers} \cite{fbc2005}.
The initial conditions are
\begin{equation}
q(y,0)\equiv-\sin(\pi y)+\hat{u}_2\sin(2\pi y).
\label{eq:bft}
\end{equation}
For the first test we choose $\nu=10^{-4}$, $\vec{k}=\uv{1}$, $\hat{u}_2=0$
and use biperiodic boundary conditions for $\x\in[0,1]^2$. In each case the problem is initialized with
$K=4\times1$ grid of a specified degree $\pode$. A BDF3/Ext3 scheme is used for the time
derivative and advective terms, respectively. We initialize only
at $t=t^0$ for this problem and integrate using a BDF1/Ext1 scheme to
provide values at $t^1$ and $t^2$.  Two cases are considered: a nonadaptive case and an
adaptive case with a maximum refinement level of $\ell_{\rm{max}}=3$. In the nonadaptive case,
the $x^1$-coordinates of the element vertices are at $x^1=0,\pm0.05,\pm1$, whereas
in the adaptive case, the elements are initially uniform. The derivative error criterion is
used in this problem, and the tolerance and coarsening multiplication are $1.0$ and $0.5$,
respectively.
Table \ref{t:burgersstat_front_na} presents the nonadaptive
results, together with the results of the nonadaptive runs from \cite{mavriplis94}. The quantity
$|\PD{u}{x}|_{\rm max}$ is the maximum of the derivative, and $T_{\rm max}$
is the time at which this occurs. We have verified that the $\Ltwo$ error of the solution
is consistent with the error in the derivative. We note that
the $\pode=5$ case is considerably worse than the 
results presented in \cite{mavriplis94}. This may be due to differences between the bases
used in the two methods \cite{basdevant86}. We note that the errors in 
$T_{\rm max}$ for our nonadaptive case are uniformly better than those in \cite{mavriplis94},
while the errors in $|\PD{u}{x}|_{\rm max}$ are comparable for $\pode>5$.
\begin{table}
\renewcommand{\arraystretch}{.7}
\begin{tabular}{c|l|l||l|l}
\hline
\hline
&\multicolumn{2}{c||}{Mavriplis}
&\multicolumn{2}{c}{\gaspar} \\
\multicolumn{1}{c|}{$\pode$}
&\multicolumn{1}{c|}{ $T_{\rm max}$}
&\multicolumn{1}{c||}{ $|\PD{u}{x}|_{\rm max}$}
&\multicolumn{1}{c|}{ $T_{\rm max}$}
&\multicolumn{1}{c}{ $|\PD{u}{x}|_{\rm max}$} \\
\hline
5  & 0.53745 & 167.227 & 0.5320  & 228.38977 \\ 
9  & 0.50611 & 154.019 & 0.51074 & 148.04258 \\ 
13 & 0.51103 & 151.496 & 0.51072 & 151.69874 \\ 
17 & 0.51071 & 152.076 & 0.51045 & 152.09104 \\ 
21 & 0.51023 & 152.004 & 0.51047 & 151.99624 \\ 
\hline 
{\rm Analytic} & 0.51047 & 152.00516 & - & - \\
\hline
\end{tabular}
\caption{Nonadaptive results from the stationary Burgers front problem. }
\label{t:burgersstat_front_na}
\bigskip
\end{table}

\begin{table}
\renewcommand{\arraystretch}{.7}
\begin{tabular}{c|l|l||l|l||l|l}
\hline
\hline
&\multicolumn{2}{c}{Adaptive}
&\multicolumn{2}{c}{Reference}
& \multicolumn{2}{c}{Control} \\
\multicolumn{1}{c|}{$\pode$}
&\multicolumn{1}{c|}{ $T_{\rm max}$}
&\multicolumn{1}{c||}{ $|\PD{u}{x}|_{\rm max}$}
&\multicolumn{1}{c|}{ $T_{\rm max}$}
&\multicolumn{1}{c}{ $|\PD{u}{x}|_{\rm max}$}
&\multicolumn{1}{c|}{ $T_{\rm max}$}
&\multicolumn{1}{c}{ $|\PD{u}{x}|_{\rm max}$} \\
\hline
5  & 0.52679 & 224.36164 &--       & --        & 0.52674 & 224.37214\\
9  & 0.51095 & 153.39634 & 0.52635 & 227.53596 & 0.51095 & 153.39633\\
13 & 0.51030 & 150.03130 & 0.51219 & 181.02024 & 0.51030 & 150.03130\\
17 & 0.51048 & 152.25110 & 0.51082 & 149.57372 & 0.51048 & 152.25110\\
21 & 0.51047 & 152.00556 & 0.51021 & 147.22940 & 0.51047 & 152.00565\\
\hline
{\rm Analytic} & 0.51047 & 152.00516 & - & - & - & -\\
\hline
\end{tabular}
\caption{Adaptive, reference, and control results from the stationary Burgers front problem. }
\label{t:burgersstat_front_ad}
\bigskip
\end{table}


In Table \ref{t:burgersstat_front_ad} we present the results from the adaptive case and the so-called reference and control solutions. The \emph{reference} solution
runs on a uniform grid with the same number of elements as the adaptive solution at $T_{\rm max}$.
Thus, it offers a solution computed with roughly as many \dof\ as the 
adaptive solution, and hence requires about the same computational effort. Clearly, 
the ability to resolve the front is critical in this case; the reference solution for $\pode=5$ 
actually diverges, and good solutions are not produced until $\pode>13$. 
The control solutions, as expected, are all nearly identical to the adaptive solutions.
This fact suggests that our refinement criteria enable the adaptive mesh to capture the formation
of the front accurately, at a significantly reduced cost.

\subsubsection{N-wave problem}
\mpar{While the 1D version of this problem can be also be considered \cite{witham74}, we choose
to use the N-wave problem in a fully 2D setting.}
The 2D Cole-Hopf transformation 
\be
\u=-2\nu\grad\ln\chi.
\label{eq:colehopf}
\ee
transforms \eq{eq:burgers} into a heat equation for $\chi$.
Choosing a source solution \cite{witham74}
$$
\chi(\x,t)=1+\frac{a}{t}\exp-\frac{(\x-\x^{\,0})^2}{4\nu t},
$$
we obtain the solution to \eq{eq:burgers} immediately from \eq{eq:colehopf}:
\be
\u(\x,t)=\frac{\x-\x^{\,0}}{t}\frac{a}{a+t\exp((\x-\x^{\,0})^2/4\nu t)}.
\label{eq:nwave}
\ee
This is a radial N-wave, where $\x^{\,0}=(\uv{1}+\uv{2})/2$ is the location of the center. This solution is singular as $t\to0$, so we
initialize at a finite time $t^0=5\times10^{-2}$. For this test problem, we choose $\nu=5\times10^{-3}$ and
$a=10^4$.  
Dirichlet boundary conditions \eqref{e:BoCo}
on $\pdomain=[0,1]^2$ are imposed at each timestep by using \eq{eq:nwave} on
$\partial\pdomain$. The initial grid consists of $K=4\times4$ elements, and we consider
only the full\mpar{What does ``full'' mean here?} adaptive case with $\ell_{\rm{max}}=4$. This $\ell_{\rm{max}}$
provides a control mesh that is
expensive to compute on; hence, we 
examine the temporal and spatial convergence of only the
adaptive solutions. The refinement criteria are the same as in \sr{sec:linadvecttest}. 

\begin{figure}\begin{center}
\mbox{\includegraphics*[width=.5\textwidth]%
{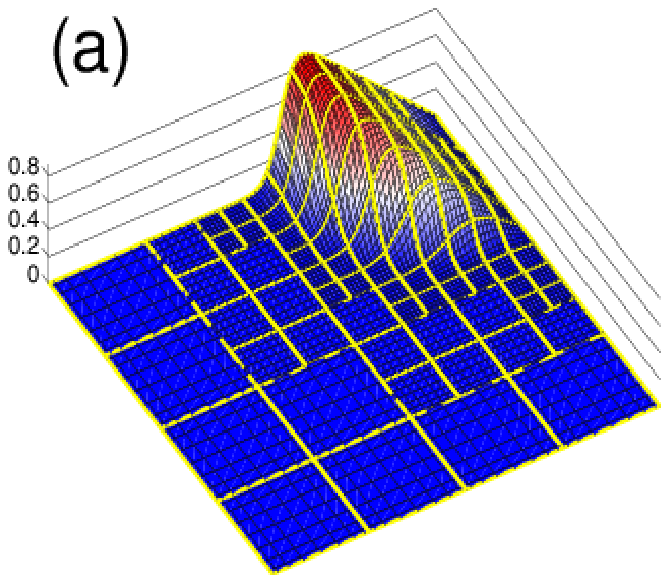}
\includegraphics*[width=.5\textwidth]%
{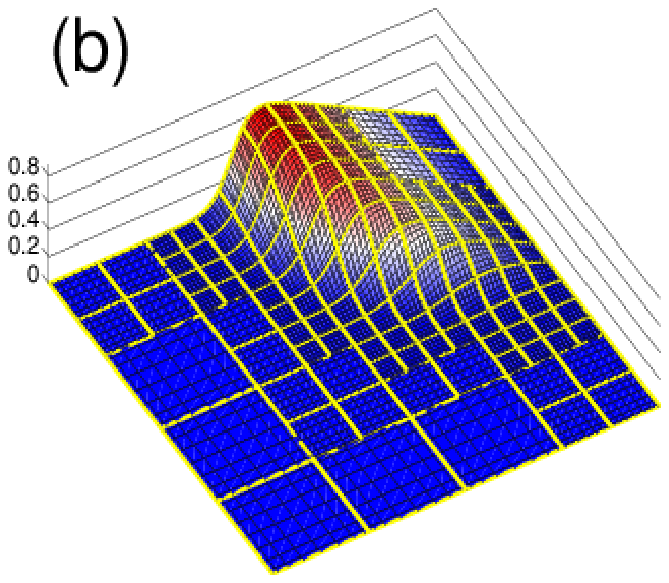}}
\mbox{\includegraphics*[width=.5\textwidth]%
{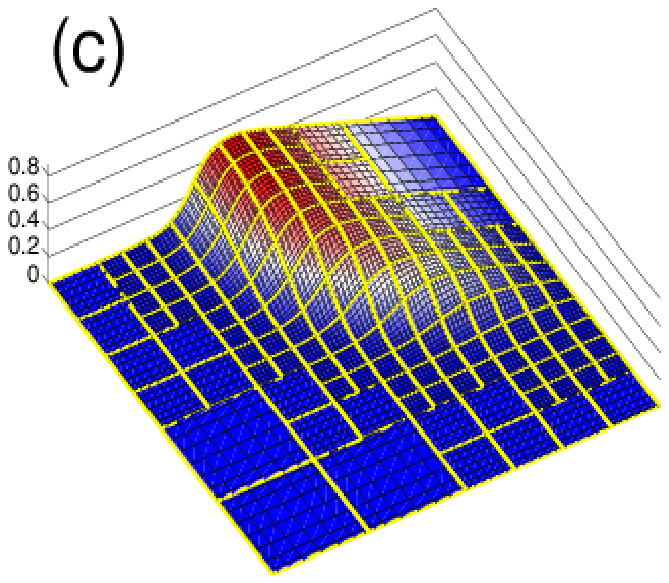}
\includegraphics*[width=.5\textwidth]%
{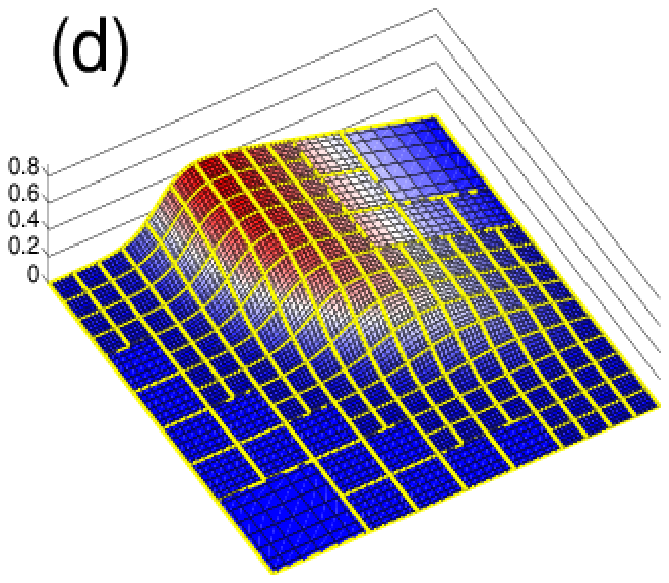}}
\mbox{\includegraphics*[width=.5\textwidth]%
{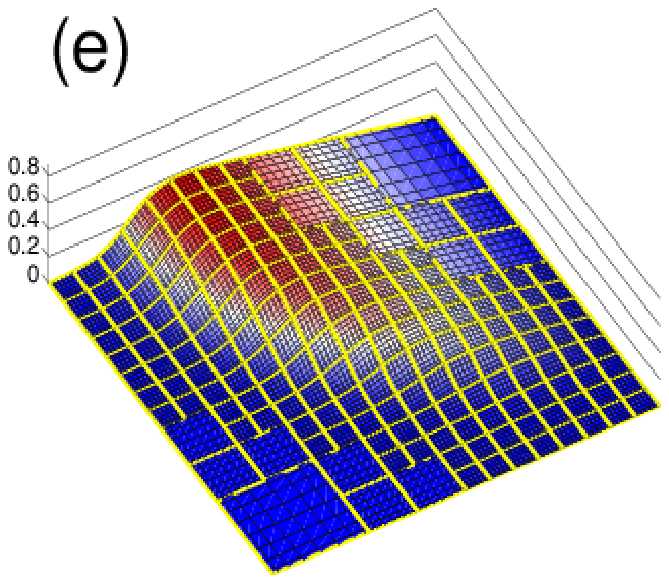}
\includegraphics*[width=.5\textwidth]%
{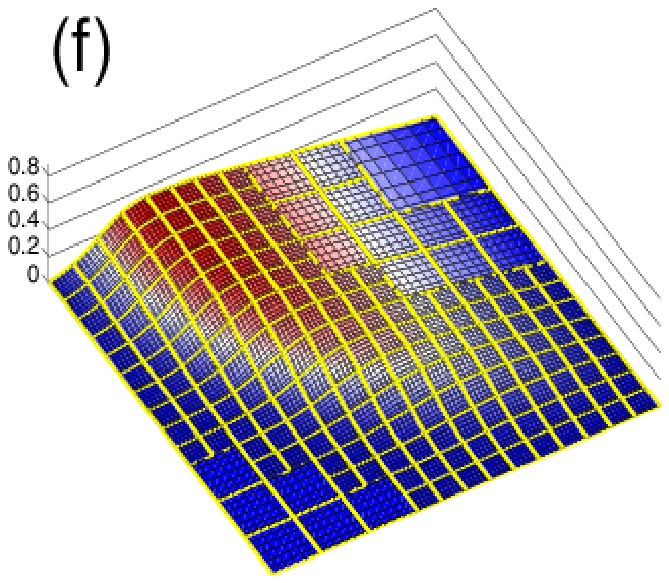}}
\end{center}
\caption{Surface plots of $u^1(\x,t)$ solved from
\eqref{eq:burgers}, showing $\x\in[\half,1]^2$ and $K/4=88$,
$121$, $139$, $172$, $181$, $190$ as $t=0.18$, $0.33$, $0.48$,
$0.65$, $0.81$, $1.00$.  Here $\ua=\u$, $\nu=5\times10^{-3}$,
$\pode=8$, and \eqref{eq:nwave} is the initial condition.}
\label{fig:nwave}
\end{figure}

\fig{fig:nwave} presents a time series of a typical N-wave solution 
illustrating the refinement patterns characteristic of all the runs.
For simplicity only one quadrant of the axially symmetric solution is presented.
As the front propagates
outward, the grid refines to track it, while in the center where the 
velocity components are planar, the grid coarsens. Note that the front does not steepen 
in this problem, as it does in the planar front problem (\sr{sec:burgerstest}); it 
simply decays as it moves outward.

\begin{figure}
\begin{center}
\includegraphics[scale=.60]{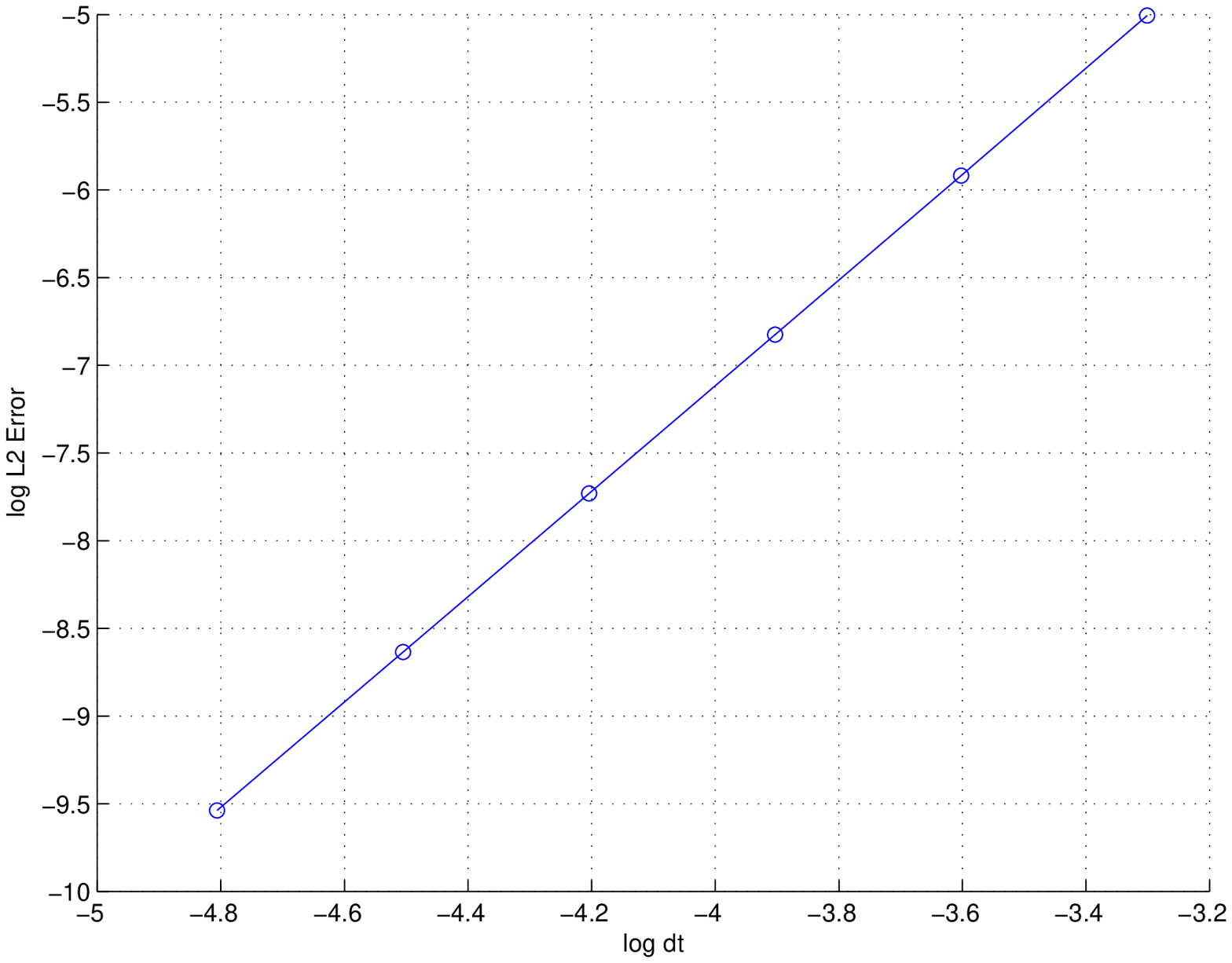}
\end{center}
\caption{Plot of $\log_{10}(\dlebnorm{\unum-\unum_{\rm{a}}}/\dlebnorm{\unum_{\rm{a}}^0})$ for the N-wave problem vs $\log_{10}\varDelta t$, for different $\pode=14$. The slope of the line is $3.01$.
\label{fig:nwave_dtconv}
}
\end{figure}

In considering the time convergence, we set $\pode=14$ and vary 
$\varDelta t$ to produce the error plot in \fig{fig:nwave_dtconv}.
As was the case previously, each point in this plot reflects a single 
run integrated at the indicated fixed 
$\varDelta t$ and the error at $t=0.11$. 
This integration time interval was chosen to provide a number
of refinement and coarsening events. 
Nevertheless, the solution converges with $\varDelta t$, 
at order (slope) $3.01$. 

\begin{figure}
\begin{center}
\includegraphics[scale=.60]{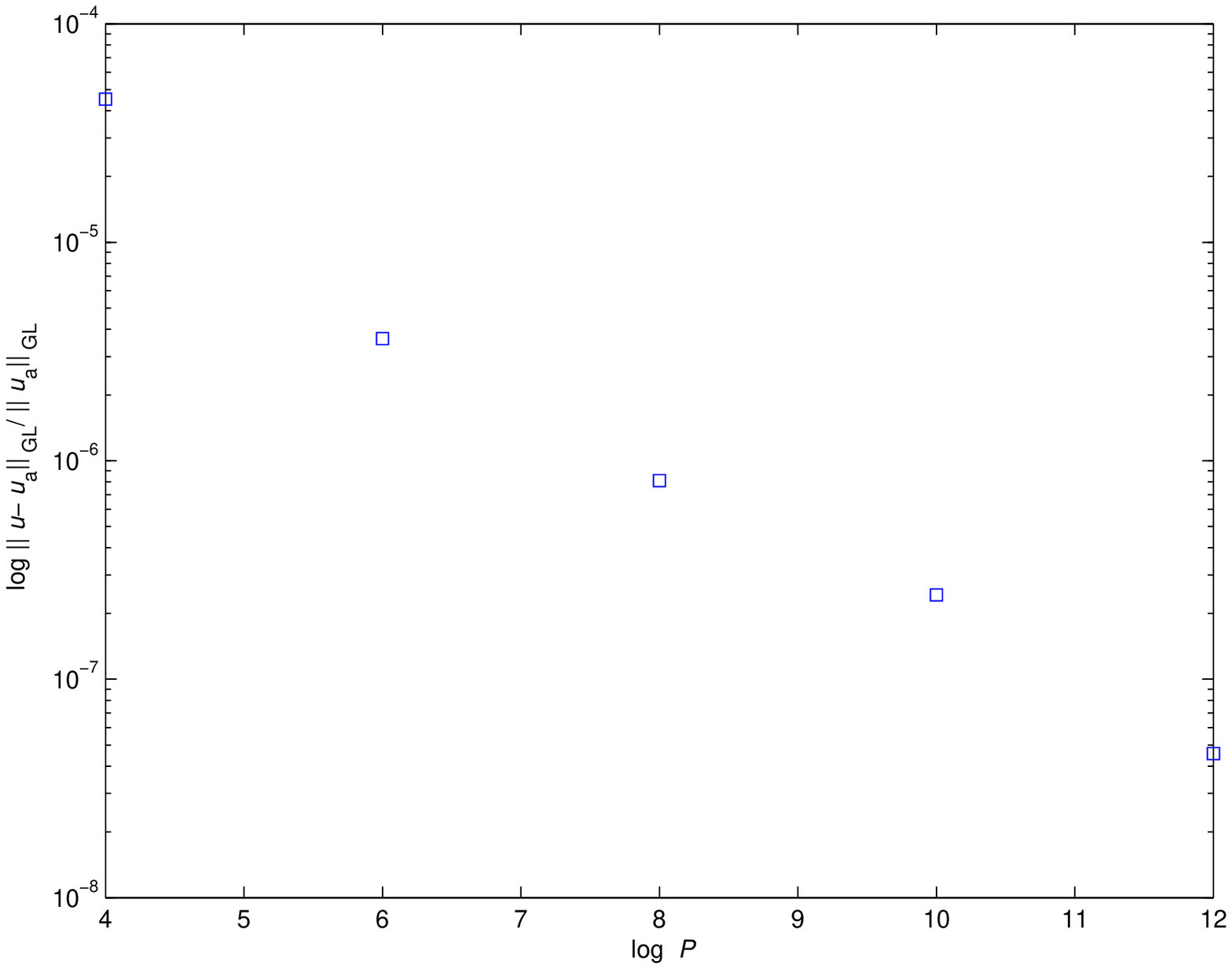}
\end{center}
\caption{Plot of $\log_{10}(\dlebnorm{\unum-\unum_{\rm{a}}}/\dlebnorm{\unum_{\rm{a}}^0})$, as a function of $\pode$ for the N-wave test. 
\label{fig:nwave_spconv}
}
\end{figure}

In order to check spatial convergence, the solution is integrated to 
$t=0.11$ by using a variable $\varDelta t$ and $\pode$ and fixed
Courant number of $0.15$. 
In \fig{fig:nwave_spconv} we present
the final $\Ltwo$ error vs $\pode$. As with the linear 
advection case, the error behaves spectrally for a finite time
integration.

\section{Discussion and conclusion}
\label{sec:conclusion}

We have presented an overview of the \gaspar code, concentrating on 
the continuous Galerkin discretization of a single 
(generalized advection-diffusion) equation to illustrate
the construction of the weak and collocation operators and to
highlight aspects of the code design. We have provided a detailed
description of the underlying mathematics and constructs
for connectivity and a new adaptive mesh refinement algorithm that 
are used in the code, and we have shown how they maintain continuity between
conforming and nonconforming elements. In the process, we have
established a consistent set of rules for refinement and coarsening.
In addition, we have described
the refinement criteria and using them, have demonstrated 
with several test problems the ability of DARe to model accurately a variety of
2D structures that evolve as a result of linear and nonlinear advective 
and dissipative dynamics.

The test problems suggest that for resolving isolated structures,
dynamic adaptive refinement can offer a substantial computational
savings over a either a pseudo-spectral or conforming spectral-element
method for the same control resolution.  But for turbulent flows,
how likely is it that only a few isolated structures will exist? And
how dependent is the time evolution of the flow on these structures,
such that, if they are resolved, the statistical properties of the overall
flow will be preserved using DARe?  One aspect of \gaspar\ that can help
answer these questions is that the fields solved for need not
be those on which adaption criteria directly operate. The
user is free to specify any functionals of the fundamental fields
(velocity, pressure, etc.) for use in tagging elements for refinement
or coarsening. For example, while the velocity is actually solved for
in \eq{eq:burgers}, the adaption criteria might operate on kinetic
energy, vorticity, or enstrophy. Arguably, some fully developed turbulent
flows viewed in terms of the fundamental fields may be too intricate
to benefit from DARe. Nevertheless, when viewed w.r.t.\
some functional, some relevant structures, when resolved,
may allow for accurate simulation of the significant dynamics of the
overall flow.

\appendix*
\section{Spectral-element formalism}
\label{app:sem}
\subsection{Canonical 1D element}

In this appendix we summarize our notation and results from the SEM literature.
Any dependent variable $u=u(\xi)$ in the domain $\xi\in[-1,1]$
may be approximated by its projection $\opr{P}_{\pode}u$ on the space
$\polynomialsset{\pode}$ of polynomials of degree $\pode$,
using $u$-values on any $\nopo\equiv\pode+1$ distinct nodal points $\Gn{j}$:
\be
u=\opr{P}_{\pode}u+\errterm{\pode}u\approx
\opr{P}_{\pode}u\equiv\sum_{j=0}^{\pode}u(\Gn{j})\phi_j,
\label{e:1dproj}
\ee
where $\errterm{\pode}u$ is the pointwise error and
\be
\phi_j(\xi)\equiv\prod_{j'\neq j}%
\frac{\xi_{\hphantom{j}}-\Gn{j'}}{\Gn{j}-\Gn{j'}}%
\tendsto{\xi}{\Gn{j''}}\delta_{j,j''},
\label{e:lagrinte}
\ee
denotes the $\nopo$ unique Lagrange interpolating polynomials
of degree $\pode$.
Now let the $\Gn{j}$ be the {\rm Gauss-Lobatto-Legendre}
quadrature nodes
\be
\Gn{j}\equiv\text{$(j+1)$th least root of
$(1-\xi^2)\frac{\d\Leg{\pode}}{\d\xi}$},
\label{e:GLLqn}
\ee
where
$\Leg{j}$ is the standard Legendre polynomial of degree $j$ and norm
$(j+\half)^{-\frac{1}{2}}$.
In this case
\be
\phi_{j'}(\xi)=\Gw{j'}\sum_{j=0}^\pode\Leg{j}(\Gn{j'})\Leg{j}(\xi)/%
\sum_{j''=0}^\pode\Gw{j''}\Leg{j}(\Gn{j''})^2,
\label{e:lagrinteex}
\ee
and the integral
\be
\mv[1]{u}\equiv\int_{-1}^1u(\xi)\d\xi=\sum_{j=0}^\pode\Gw{j}u(\Gn{j})+%
\opr{R}_{\pode}u(\xi'),
\label{e:1dquad}
\ee
where
\be
\Gw{j}\equiv2/\pode\nopo\Leg{\pode}(\Gn{j})^2
\label{e:GLLqw1}
\ee
denotes the Gauss-Lobatto-Legendre weights,
$\opr{R}_\pode\equiv%
-2^{2\pode+1}\frac{\pode^3\nopo(\pode-1)!^4}{(2\pode+1)(2\pode)!^3}(\d/\d\xi)^{2\pode}$
is the residual operator \cite{WeisLQ} and $\xi'\in\ropen-1,1\lclose$.
Then the mean-square error is bounded as
\be
\mv[1]{(\errterm{\pode}u)^2}\propto\pode^{1-2Q}\sum_{q=0}^Q\|u^{(q)}\|^2
\label{e:errorbound}
\ee
for any order $Q$
of square-integrable derivative.
Thus if $u$ is infinitely smooth then $\opr{P}_{\pode}u$
converges to $u$ \emph{spectrally}.

\subsection{General 1D spectral elements}

Now let $[-1,1]$ be subdivided and covered by $K^1$ disjoint 1D elements
$\ropen\pn{k-1},\pn{k}\lclose\equiv\elbox{k}^1\equiv\comap{k}(\ropen-1,1\lclose)$,
where
\be
\comap{k}(\xi)\equiv\pn{k-1}+\half\elsiz[1]{k}(1+\xi),\qquad k\isin{1}{K^1},
\label{e:comap}
\ee
is a coordinate map with inverse
$\comap{k}^{-1}(x)=2(x-\pn{k-1})/\elsiz[1]{k}-1$.
(Other invertible maps may sometimes be preferable.)
For $\pn{0}\equiv-1$, $\pn{K^1}\equiv1$ and positive
$\elsiz[1]{k}\equiv\pn{k}-\pn{k-1}$ one has
$\elbox{k}^1\bigcap\elbox{k'}^1=\varnothing$ if $k\neq k'$ and
\be
[-1,1]=\union[K^1]{k=1}\elbox[\bar]{k}^1.
\label{domdecom1}
\ee
Then $u$ may be approximated by its projections $\opr{P}_{k,\pode}u$ on the space
$\pwisepolysset{\lavec{\eld}^1}{\pode}$ of piecewise polynomials of degree $\pode$ on the
$\elbox{k}^1$,
by using $u$-values on the $K^1\nopo$ mapped nodal points
\be
\pn{j,k}\equiv\comap{k}(\Gn{j})
\label{e:npmapped1}
\ee
generalized from \eqref{e:GLLqn}.
That is, \eqref{e:1dproj} generalizes to
\be
u=\sum_{k=1}^{K^1}\left(\opr{P}_{k,\pode}u+\errterm{k,\pode}u\right),\qquad
\opr{P}_{k,\pode}u\equiv\sum_{j=0}^{\pode}u(\pn{j,k})\phi_{j,k},
\label{e:1dprojmapped}
\ee
where
$\errterm{k,\pode}u\equiv\errterm{\pode}(u\circ\comap{k})\circ\comap{k}^{-1}$, and
\eqref{e:lagrinte} generalizes to
\be
\phi_{j,k}(x)\equiv1_{\elbox[\bar]{k}^1}(x)\phi_j\circ\comap{k}^{-1}(x)%
\tendsto{x}{\pn{j',k'}}\begin{cases}1,&x_{j,k}=x_{j',k'},\\0,&\text{otherwise}.\end{cases}
\label{e:limapped1}
\ee
We have adopted the additional notations
$f\circ g(x)\equiv f(g(x))$
and $1_{\set{S}}(x)\equiv1$ $(x\in\set{S})$, $0$ (otherwise).
Then \eqref{e:1dquad} generalizes to
\be
\mv[1]{u}=\sum_{k=1}^{K^1}\int_{\pn{k-1}}^{\pn{k}}u(x)\d x,\quad
\int_{\pn{k-1}}^{\pn{k}}u(x)\d x=%
\sum_{j=0}^\pode\Gw{j,k}u(\pn{j,k})+\opr{R}_{k,\pode}u(\pn{k}'),
\label{e:1dequad}
\ee
where \eqref{e:GLLqw1} generalizes to
\be
\Gw{j,k}\equiv|\frac{\d}{\d\xi}\comap{k}(\Gn{j})|\Gw{j}
,\label{e:GLLqwmapped}
\ee
$\opr{R}_{k,\pode}u\equiv(\elsiz[1]{k}/2)^{2\pode+1}\opr{R}_{\pode}(u\circ\comap{k})\circ\comap{k}^{-1}$
and $\pn{k}'\in\elbox{k}^1$.

\subsection{General $d$-dimensional spectral elements}

Generalizing \eqref{domdecom1},
assume the $d$-dimensional problem domain $\pdomain$ can be partitioned into
$K$ disjoint images of $[-1,1]^d\equiv\setdef{\vec{\xi}}{\xi^\si\in[-1,1]}$ as
\be
\bar{\pdomain}=\union[K]{k=1}\elbox[\bar]{k},
\label{e:domdecom}
\ee
where $\elbox{k}\equiv\comap[\vec]{k}(\ropen-1,1\lclose^d)$ has diameter
\be
\elsiz{k}\equiv\max_\si\max_{\x,\x'\in\elbox[\bar]{k}}|x^\si-{x'}^\si|
\label{e:elsiz}
\ee
and
$\comap[\vec]{k}(\vec{\xi})$ has inverse $\comap[\vec]{k}^{-1}(\x)$ but is not
necessarily linear.
Now generalizing \eqref{e:1dprojmapped}, one may approximate a field
$u(\x)$ by its projections
$\opr{P}_{k,\vec{\pode}\,}u$ on the space
$\pwisepolysset{\lavec{\eld}}{\vec\pode}$ of piecewise polynomials of degree $\pode^\si$
in coordinate $x^\si$ on the $\elbox{k}$,
using $u$-values on the $\nopo[K,\vec{\pode}]\equiv K\prod_{\si=1}^d\nopo[\pode^\si]$ mapped nodes
\be
\pn[\vec]{\vec{\jmath},k}\equiv\comap[\vec]{k}(\Gn[\vec]{\vec{\jmath}}),
\label{e:npmapped}
\ee
generalized from \eqref{e:npmapped1}, where
$\Gn{\vec{\jmath}}^\si\equiv\Gn{\jmath^\si}$ are $d$-dimensional GLL nodes.
That is, \eqref{e:1dprojmapped} generalizes to
\be
u\approx\opr{P}_{\lavec{\eld},\vec\pode\,}u\equiv%
\sum_{k=1}^K\opr{P}_{k,\vec{\pode}\,}u,\qquad
\opr{P}_{k,\vec{\pode}\,}u\equiv\sum_{\vec{\jmath}\in\set{J}}%
u(\pn[\vec]{\vec{\jmath},k})\phi_{\vec{\jmath},k},
\label{projmapped}
\ee
where
$\set{J}\equiv\setdef{\vec{\jmath}\,}{\jmath^\si\isin{0}{\pode^\si}}$,
\eqref{e:limapped1} generalizes to
\be
\phi_{\vec{\jmath},k}(\x)\equiv%
1_{\elbox{k}}(\x)\phi_{\vec{\jmath}}\circ\comap[\vec]{k}^{-1}(\x)%
\tendsto{\x}{{\pn[\vec]{\vec{\jmath}',k'}}}\begin{cases}%
1,&\pn[\vec]{\vec{\jmath},k}=\pn[\vec]{\vec{\jmath}',k'}\\
0,&\text{otherwise,}\end{cases}
\label{e:limapped}
\ee
and
$\phi_{\vec{\jmath}}(\vec{\xi})\equiv%
\prod_{\si=1}^d\phi_{\jmath^\si}(\xi^\si)$.
The appropriate approximation of a vector
$$
\u=\sum_{\si=1}^d u^\si\uv{\si}\approx\opr{P}_{\lavec{\eld},\vec\pode\,}\u%
=\tr{\smash{\vec{\lavec{\phi}}}}\unum
$$ uses $\vec{\lavec{\phi}}$ with entries
\be
\vec{\phi}_{\vec{\jmath},k}^\si\equiv\phi_{\vec{\jmath},k}\uv{\si}
\label{e:limappedv}
\ee
and $\unum$ with entries
$u_{\vec{\jmath},k}^\si\equiv u^\si(\pn[\vec]{\vec{\jmath},k})$.
For scalars $u$ \eqref{e:1dequad} generalizes to
\be
\begin{split}
\mv{u}&\equiv\idotsint\limits_\pdomain u(\x)\d^d\x%
=\sum_{k=1}^K\idotsint\limits_{\elbox{k}}u(\x)\d^d\x\\
&\approx\sum_{k=1}^K\sum_{\vec{\jmath}\in\set{J}}%
\Gw{\vec{\jmath},k}u(\pn[\vec]{\vec{\jmath},k})\equiv\mv{u}_\GL,
\end{split}
\label{e:quad}
\ee
where \eqref{e:GLLqwmapped} generalizes to
\be
\Gw{\vec{\jmath},k}\equiv|\det\grad_{\vec{\xi}\,}\comap[\vec]{k}(\Gn[\vec]{\vec{\jmath}})|\prod_{\si=1}^d\Gw{\jmath^\si}.
\label{e:GLLqw}
\ee
Finally, variational formulation depends on the inner product from \eqref{e:quad}:
\be
\ipc{u}{v}\equiv\mv{uv}\approx\sum_{k=1}^K\sum_{\vec{\jmath}\in\set{J}}%
\Gw{\vec{\jmath},k}u(\pn[\vec]{\vec{\jmath},k})v(\pn[\vec]{\vec{\jmath},k})\equiv%
\ipc{u}{v}_\GL
\label{e:innpro}
\ee
for scalars, $\ipc{\u}{\v}\equiv\sum_{\si=1}^d\ipc{u^\si}{v^\si}$ for vectors,
$\ipc{\Vec{\Vec{u}}}{\Vec{\Vec{v}}}\equiv%
\sum_{\si,\si'=1}^d\ipc{u^{\si,\si'}}{v^{\si,\si'}}$ for tensors, and so forth.
%

\section*{Acknowledgments}

We thank Rich Loft, Peter Sullivan, Steve Thomas and Joe Tribbia
for help at the onset of this work in using spectral element methods,
and Huiyu Feng and Catherine Mavriplis for several useful discussions.
Computer time was provided by NSF under sponsorship of the National
Center for Atmospheric Research. The third author was supported by the U.S.
Department of Energy under Contract W-31-109-ENG-38.

%

%
\end{document}